\documentclass[twoside]{article}
\usepackage[margin=1in]{geometry}
\usepackage[english]{babel}
\usepackage{amsmath}
\usepackage{amsfonts}
\usepackage{amssymb}
\usepackage{amsthm}
\usepackage[nocompress]{cite}
\usepackage{graphicx}
\usepackage{epstopdf}
\usepackage{algorithm}
\usepackage{algorithmic}
\usepackage{footnote}
\usepackage{multirow}
\usepackage{hyperref}
\usepackage[capitalise]{cleveref}
\usepackage{bm}
\usepackage{bbm}

\crefformat{equation}{(#2#1#3)}
\newtheorem{theorem}{Theorem}[section]
\newtheorem{remark}[theorem]{Remark}

\crefname{remark}{Remark}{Remarks}
\crefname{assumption}{Assumption}{Assumptions}

\makeatletter
\renewcommand\paragraph{%
  \@startsection{paragraph}
    {4}
    {\z@}
    {3.25ex \@plus1ex \@minus.2ex}
    {-1em}
    {\normalfont\normalsize\bfseries\maybe@addperiod}%
}
\newcommand{\maybe@addperiod}[1]{%
  #1\@addpunct{.}%
}
\makeatother

\usepackage{fancyhdr}
\title{Geometry-based approximation of waves in complex domains\thanks{This research was funded in part by the Austrian Science Fund (FWF) projects 10.55776/F65 (MN, IP), 10.55776/P33477 (MN, IP), and 10.55776/ESP519 (MN).}}
\newcommand{\email}[1]{\protect\href{mailto:#1}{#1}}
\author{Davide Pradovera\thanks{Department of Mathematics, KTH Royal Institute of Technology, Lindstedtsv{\"a}gen~25, 11428 Stockholm, Sweden (\email{davidepr@kth.se}).}\ \thanks{Faculty of Mathematics, University of Vienna, Oskar-Morgenstern-Platz~1, 1090 Vienna, Austria
  (\email{monica.nonino@univie.ac.at}, \email{ilaria.perugia@univie.ac.at}).}
\and Monica Nonino\footnotemark[3]
\and Ilaria Perugia\footnotemark[3]}

\newcommand{\uapp}{\widetilde{u}}
\newcommand{\ufem}{u_{\textup{FE}}}

\newcommand{\bR}{\mathbb R}

\newcommand{\va}{\bm a}

\newcommand{\vx}{\bm x}
\newcommand{\vxi}{\bm \xi}
\newcommand{\vy}{\bm y}
\newcommand{\vz}{\bm z}
\newcommand{\vzero}{\bm 0}
\newcommand{\ks}{\mu}

\DeclareMathOperator*{\argmin}{arg\,min}
\newcommand{\abs}[1]{\left|#1\right|}
\newcommand{\norm}[1]{\left\|#1\right\|}
\newcommand{\one}{\mathbbm 1}

\ifpdf\hypersetup{pdftitle={Geometry-based approximation of waves in complex domains},pdfauthor={D. Pradovera, M. Nonino, and I. Perugia}}\fi

\begin{document}

\maketitle

\begin{abstract}
We consider wave propagation problems over 2-dimensional domains with piecewise-linear boundaries, possibly including scatterers. We assume that the wave speed is constant, and that the initial conditions and forcing terms are radially symmetric and compactly supported. We propose an approximation of the propagating wave as the sum of some special space-time functions. Each term in this sum identifies a particular field component, modeling the result of a single reflection or diffraction effect. We describe an algorithm for identifying such  components automatically, based on the domain geometry. To showcase our proposed method, we present several numerical examples, such as waves scattering off wedges and waves propagating through a room in presence of obstacles. Software implementing our numerical algorithm is made available as open-source code.
\end{abstract}



\section{Introduction}
\label{sec:intro}
The discretization of numerical models for the simulation of complex phenomena results in high-dimensional systems to be solved, usually at an extremely high cost in terms of computational time and storage memory. Among these models, wave propagation problems represent an extremely interesting topic: relevant applications can be found, e.g., in the field of array imaging, where acoustic, electromagnetic, and elastic waves in scattering media are modeled by the \emph{reflectivity} coefficient, which is often unknown. Some examples in this direction can be found in \cite{Borcea3,Borcea,Borcea2,Tournier2019}, where inverse scattering problems are used to infer the reflectivity of one or more scatterers embedded either in a known and smooth medium, or in a randomly inhomogeneous medium. Another example of application of wave propagation problems is numerical acoustics, where the goal is to simulate the propagation of sound in a room, in presence of obstacles and walls with different absorbing and/or reflecting properties. See, e.g., \cite{potter2022numerical,svensson_savioja_overview} for some frequency-domain examples of this.

Our focus here are problems in the \emph{time domain}, whose numerical simulation is expensive, mainly because one needs to use both a fine spatial mesh and a carefully chosen time step in order to satisfy the CFL condition \cite{Cohen2008,HairerLubich}. In the interest of making these simulations feasible, model order reduction (MOR) \cite{benner2017model,glas2020reduced,hesthaven2021structure} represents a promising framework, whose goal is to reduce the computational cost of solving the problem of interest.

In this context, it is well known \cite{greif2019decay} that wave propagation problems are characterized by a slowly decaying Kolmogorov $n$-width. Because of this, classical linear-subspace MOR methods are not able to reproduce the behavior of the wave propagation without relying on a very high-dimensional linear manifold. This makes linear surrogate models unappealing, since they do not yield significant speed-ups. In recent years, methods that rely on nonlinear and/or hybrid space-time approaches have been proposed. See, e.g., \cite{Cagniart2019,kleikamp_nonlinear_2022,reiss2018shifted} and references therein. 

In this work, we focus on wave propagation over 2-dimensional  spatial domains, possibly including obstacles. We aim at proposing a new approximation framework, which satisfies the main goal of MOR techniques, that is, making the computational cost of the numerical simulations more feasible.
We limit our investigation to domains with piecewise-linear boundaries and a constant wave speed. The initial conditions and forcing terms are assumed to be compactly supported and radially symmetric around a ``source point''. Under these assumptions, and following ideas similar to those in \cite{potter2022numerical,svensson_savioja_overview}, we propose to construct a surrogate model for the solution of the target problem as a superposition of some special space-time terms, which we call ``field components''.

Each field component models a reflection or diffraction effect, and is characterized by:
\begin{itemize}
	\item a space-time propagation term, related to the \emph{free-space radially symmetric solution} of the wave equation;
	\item a spatial indicator function, determining the \emph{spatial support} of each component;
	\item a nonlinear spatial term encoding the \emph{angular modulation} of the component, which is crucial when modeling diffraction effects.
\end{itemize}
To compute the ``angular modulation'' term, we leverage concepts from the frequency-domain \emph{geometric theory of diffraction} \cite{keller1962geometrical,kouyoumjian1974uniform}. However, in order to obtain a time-domain diffraction model that is suitable for our purposes, we first need to adapt the available results, effectively developing a novel time-domain diffraction model in the process.

The number of field components appearing in the surrogate is determined by the number of reflection and diffraction effects that are required to faithfully approximate the target wave, which ultimately depends on the geometry of the physical domain.

Among the advantages of the proposed approach, we mention the fact that each field component is separable into time-radial and angular components (in the ``polar coordinates'' sense). As we will see, we can leverage this to drastically reduce the computational cost and storage memory requirements of our approximation strategy.

The rest of the paper is structured as follows. In \cref{sec:problem} we present the problem of interest. In \cref{sec:approx} we introduce the main ingredients of our method, and we present our surrogate modeling algorithm. In \cref{sec:reflect,sec:diffract} we detail how we model reflection and diffraction effects, respectively. In \cref{sec:numerical} we present some numerical results to showcase our method. Both simple benchmarks (wedges) and more complicated tests (2D room model with scatterers) are considered. Some final considerations follow in \cref{sec:conclusions}.

\subsection{Target problem}\label{sec:problem}
We are interested in the numerical approximation of the solution of the wave equation in complex domains. In this work, we consider 2-dimensional domains only. However, most of our discussion generalizes to 3D. We defer a discussion on this till \cref{sec:conclusions}.

We denote by $\Omega\subset\bR^2$ the physical domain in which the wave equation is considered. We assume that $\Omega$ is either a closed polygon or a set-subtraction of polygons (to allow for multiply connected domains). We denote by $n_e$ and $n_v$ the number of edges and vertices of $\partial\Omega$, respectively. We study the propagation of waves in $\Omega$ over a given time interval of interest $[0,T]$. The model problem is the wave equation with constant (unit) wave speed:
\begin{equation}\label{eq:wave}
	\begin{cases}
		\partial_{tt}u(\vx,t)=\Delta u(\vx,t)+f(\vx,t)\quad&\text{for }(\vx,t)\in\Omega\times(0,T),\\
		u(\vx,0)=u_0(\vx),\partial_tu(\vx,0)=u_1(\vx)\quad&\text{for }\vx\in\Omega,\\
		\partial_\nu u(\vx,t)=0\quad&\text{for }(\vx,t)\in\partial\Omega\times(0,T],
	\end{cases}
\end{equation}
with $\Delta$ the Laplacian operator, defined, in 2 dimensions, as $\Delta=\sum_{j=1}^2\partial_{x_jx_j}$. The homogeneous Neumann condition (i.e., the last equation above) models the whole boundary $\partial\Omega$ as \emph{sound-hard} \cite{Cohen2008}. More generally, all or parts of $\partial\Omega$ may be modeled as \emph{sound-soft} via a Dirichlet-type condition: $u(\vx,t)=0$. Depending on the target application, \emph{impedance} boundary conditions may also be appropriate: $Z\partial_\nu u(\vx,t)+\partial_tu(\vx,t)=0$, with $Z\geq0$.

We assume that the initial conditions $u_0$ and $u_1$, as well as the forcing term $f$, have radial symmetry around a given point. Without loss of generality, we will take such point to be the origin of $\bR^2$:
\begin{equation}
	u_0(\vx)=\eta_0(\norm{\vx}),\; u_1(\vx)=\eta_1(\norm{\vx}),\; f(\vx,t)=\eta_2(\norm{\vx},t),
\end{equation}
for all $(\vx,t)\in\Omega\times(0,T)$, with $\norm{\vx}^2=\sum_{j=1}^2 x_j^2$. We further assume that the functions $\eta_j$ have compact support, namely, that there exist $R>0$ such that $\eta_j(\rho)=0$ for all $\rho>R$ and $j=0,1,2$. Moreover, to avoid incompatibilities with the boundary conditions, for simplicity we will only consider the situation where the supports of the functions $\eta_j$ are fully contained in $\Omega$.

\section{Approximation framework}\label{sec:approx}
Before we can model boundary effects (reflection and diffraction), we need to understand how the solution $u$ would behave if no boundary were present. To this aim, we consider the wave equation in free space
\begin{equation}
	\begin{cases}
		\partial_{tt}U(\vx,t)=\Delta U(\vx,t)+f(\vx,t)\quad&\text{for }(\vx,t)\in\bR^2\times(0,\infty),\\
		U(\vx,0)=u_0(\vx)\quad&\text{for }\vx\in\bR^2,\\
		\partial_tU(\vx,0)=u_1(\vx)\quad&\text{for }\vx\in\bR^2,
	\end{cases}
\end{equation}
which we have obtained from \cref{eq:wave} by replacing $\Omega$ with the whole plane.

Due to radial symmetry (of the initial conditions and of the forcing term), we can recast the problem in polar coordinates. This allows us to define the free-space solution $\Psi$ in the radial coordinate, as the solution of
\begin{equation}\label{eq:Psi}
	\begin{cases}
		\partial_{tt}\Psi(\rho,t)=\widehat{\Delta}\Psi(\rho,t)+\eta_2(\rho,t)\quad&\text{for }(\rho,t)\in(0,\infty)\times(0,\infty),\\
		\Psi(\rho,0)=\eta_0(\rho)\quad&\text{for }\rho\in[0,\infty),\\
		\partial_t\Psi(\rho,0)=\eta_1(\rho)\quad&\text{for }\rho\in[0,\infty),
	\end{cases}
\end{equation}
where $\widehat{\Delta}$ is the Laplace operator in polar coordinates (under radial symmetry), i.e., $\widehat{\Delta}=\partial_{\rho\rho}+\frac1\rho\partial_{\rho}$, and $U(\vx,t)=\Psi(\norm{\vx},t)$ for all $\vx\in\bR^2$. For closure, it may also be necessary to provide a condition at $\rho=0$, e.g., $\partial_\rho\Psi(0,\cdot)=0$ if the data $\eta_j$ are sufficiently smooth \cite{bcRadial}. Note that, by the compact support of the initial conditions and of the forcing term, and by the finite (unit) speed of propagation of the wave equation, we have $\Psi(\rho,t)=0$ whenever $\rho>t+R$.

\begin{remark}
Generally, the free-space solution $\Psi$ is not available analytically, except for very simple choices of initial conditions and forcing term. Accordingly, in most applications, the function $\Psi$ will need to be replaced with a suitable approximation. To this effect, one can discretize \cref{eq:Psi}, e.g., with a finite element approximation (in space) and some timestepping scheme (in time). See \cref{sec:numerical} for more details on this.
\end{remark}

Our goal is to approximate, for all $(\vx,t)\in\Omega\times[0,T]$, the solution $u(\vx, t)$ of the wave equation problem in \cref{eq:wave} with a surrogate, defined as a sum of special functions $\uapp_n$ (\emph{field components}):
\begin{equation}\label{eq:sum}
u(\vx,t)\approx\uapp(\vx,t)=\sum_{n=1}^N\underbrace{\Psi(\norm{\vx-\vxi_n}+r_n,t)\one_{\Omega_n}(\vx)\zeta_n(\vx-\vxi_n)}_{\uapp_n(\vx,t)}.
\end{equation}
Therein, $\Psi$ is the above-mentioned free-space radially symmetric solution of \cref{eq:Psi}, and $\one_A$ denotes the indicator function with support $A$, i.e.,
\begin{equation}
\one_A(y)=\begin{cases}
	1\quad&\text{if }y\in A,\\
	0\quad&\text{if }y\notin A.
\end{cases}
\end{equation}
Moreover, in \cref{eq:sum}, we have introduced the following quantities:
\begin{itemize}
	\item $N$ is the number of field components used in the approximation.
	\item $\vxi_n$ is the location of the point source from which $\uapp_n$ originates.
	\item $r_n\geq 0$ is a spatial delay, which is used for the synchronization of diffraction effects.
	\item $\Omega_n\subset\Omega$ is the spatial support of $\uapp_n$.
	\item $\zeta_n:\bR^2\setminus\{\vzero\}\to\bR$ is a weight function describing the angular modulation. We assume $\zeta_n$ to be a positive-homogeneous function, i.e., $\zeta_n(\vy)=\zeta_n(\vy/\norm{\vy})$ for all $\vy\in\bR^2\setminus\{\vzero\}$.
\end{itemize}

\begin{remark}\label{rem:homogeneous}
    We refer to $\zeta_n$ as ``angular modulation'' since, in 2D, the positive-homogeneity assumption is equivalent to requiring $\zeta_n(\vy)$ to depend only on the angle between the vector $\vy$ and some reference direction, e.g., the positive $x_1$-axis.
\end{remark}

Note that, due to the finite speed of propagation of the free-space solution $\Psi$, we have that a generic term $\uapp_n(\vx,t)$ is zero whenever $t<\norm{\vx-\vxi_n}+r_n-R$, i.e., for $t$ small enough, depending on $\vx$.

The number of field components $N$ in \cref{eq:sum} will be determined based on how many boundary effects (reflections and diffractions) need to be included in $\uapp$ in order to have a good approximation of the target wave $u$. We describe a strategy for automatically identifying a good $N$ in the next section.

\subsection{Building the low-rank skeleton}\label{sec:skeleton}
Recalling that $u$ solves the wave equation in \cref{eq:wave} in the domain $\Omega$, we use the first term in \cref{eq:sum}, namely, $\uapp_1$, to approximate the outgoing component of $u$, \emph{ignoring any effect due to the boundary $\partial\Omega$, except for shadows}. Then, given such $\uapp_1$, we use the other terms $\uapp_2,\ldots,\uapp_N$ to correct this first approximation. Each extra term models a single effect due to a certain portion of the boundary, specifically, an edge (reflection off that edge) or a vertex (diffraction about that vertex).

Going back to the first field component $\uapp_1$, let us define it, by providing its ``ingredients'' $\vxi_1$, $r_1$, $\Omega_1$, and $\zeta_1$, cf.~\cref{eq:sum}. We set $\vxi_1=\bm{0}$, the center of the initial condition, as well as $r_1=0$, since no delay is necessary for this first term. Then, leveraging symmetry, we set $\zeta_1\equiv 1$, which corresponds to the (physical) assumption that the propagation of $\uapp_1$ is purely radial. Finally, we set $\Omega_1$ 
(the spatial support of the first field component around $\bm{0}$) 
as the set of points that can be reached from $\bm{0}$ via a straight line without going outside $\partial\Omega$, i.e.,
\begin{equation}\label{eq:cone}
	\Omega_1=\left\{\vx\in\Omega\ :\ \tau\vx\in\Omega\ \ \forall0\leq\tau\leq 1\right\}.
\end{equation}
In summary, the first term of $\uapp$ is
\begin{equation}\label{eq:first}
	\uapp_1(\vx,t)=\Psi(\norm{\vx},t)\one_{\Omega_1}(\vx).
\end{equation}

\begin{remark}\label{rem:cone}
    From a practical viewpoint, the implementation of spatial supports like $\Omega_1$ is not trivial. In our code, we verify that an arbitrary point $\vx\in\Omega$ belongs to $\Omega_1$ by checking if there are no intersections between the line segment from $\vx$ to $\vxi_1$ and any of the edges composing $\partial\Omega$. In some sense, this may be considered akin to ``ray-tracing'' \cite{svensson_savioja_overview}. A similar approach, with minor modifications, works for the other spatial supports $\Omega_n$ related to reflected or diffracted waves, cf.~\cref{eq:reflcone,eq:diffcone} below.
\end{remark}

Then we can move to the subsequent terms $\uapp_n$, $n\geq 2$. Their expressions depend on our choice of reflection and diffraction modeling, and will be provided in the upcoming sections. Instead, in the rest of the present section we focus on understanding how large $N$ should be, in order for $\uapp$ to provide a faithful approximation of $u$. Equivalently, we want to count the number of times the wave gets reflected or diffracted at the boundary $\partial\Omega$. This is done incrementally, starting from the initial value $N=1$ (no boundary effects) and then updating this guess as more and more boundary effects get ``discovered''.

\begin{figure}[tbp]
	\centering
	\includegraphics{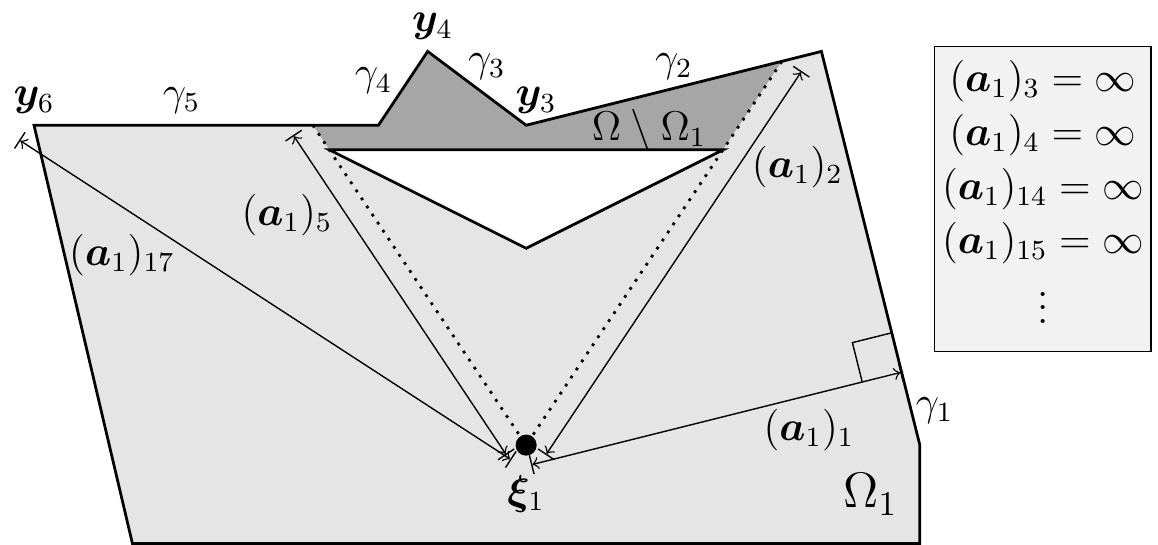}
	\caption{Computation of some timetable entries. The boundary $\partial\Omega$ has 11 sides, so that, e.g., $(\va_1)_{14}$ is related to $\vy_3$ and $(\va_1)_{17}$ is related to $\vy_6$. The shadowed area $\Omega\setminus\Omega_1$ is in darker grey.}
	\label{fig:timetable}
\end{figure}

To help us in this endeavor, we employ what we call a \emph{timetable}, which, in this work, is simply a list of vectors, each with size $n_e+n_v$. (Note that similar ideas can can also be found in \cite{lee_efficient_1988,potter2022numerical,svensson_savioja_overview} and references therein.) The timetable is built incrementally starting from an empty list, appending one new vector every time a new term is added in the sum in \cref{eq:sum}, starting from $\uapp_1$. The entries of the $n$-th timetable vector are the waiting times before $\uapp_n$ comes in contact with an edge or a vertex of $\partial\Omega$. If it is impossible for $\uapp_n$ to ``cast light'' (along a straight path) onto a certain edge or vertex (e.g., $\gamma_3$ and $\vy_4$ in \cref{fig:timetable}), then the corresponding entry in the timetable is set to $\infty$. After this, it suffices to look for the smallest not-yet-explored entry of the timetable to identify what the next term of the approximation $\uapp$ should be\footnote{Note that ties are possible when selecting the smallest timetable entry, e.g., if $\uapp_n$ reaches an edge and an adjacent vertex at the same time (for an example, see the bottom vertex and the two adjacent sides of the triangle in \cref{fig:timetable}). In such cases, each of such timetable entries must be explored \emph{in arbitrary order}, one after the other.}. Once the entry in the timetable has been explored, its value is set to $\infty$.

We start by describing how the first vector $\va_1\in\bR^{n_e+n_v}$ of the timetable (corresponding to $\uapp_1$) is computed, and how $\va_1$ allows us to identify the (geometric) features of $\uapp_2$. The vector $\va_1$ can be partitioned into edges-related part (the first $n_e$ entries) and vertices-related part (the last $n_v$ entries).

\begin{itemize}
	\item \textbf{Edge times.} Given a generic edge $\gamma_j\subset\partial\Omega$ ($j=1,\ldots,n_e$) belonging to the domain boundary, we define the corresponding entry of $\va_1$ as
	\begin{equation}\label{eq:timeE}
		\left(\va_1\right)_j=\begin{cases}
			r_1+\inf\left\{\norm{\vx-\vxi_1}:\vx\in \gamma_j\cap\overline{\Omega}_1\right\}&\text{if the set is non-empty},\\
			\infty&\text{otherwise}.
		\end{cases}
	\end{equation}
	Note that we have taken the shortest path from $\vxi_1$ to $\gamma_j$, and that we have denoted the closure of $\Omega_1$ as $\overline{\Omega}_1$.
	\item \textbf{Vertex times.} Given a generic vertex $\vy_j\subset\partial\Omega$ ($j=1,\ldots,n_v$) of the domain boundary, we set
	\begin{equation}\label{eq:timeV}
		\left(\va_1\right)_{n_e+j}=\begin{cases}
			r_1+\norm{\vy_j-\vxi_1}\quad&\text{if }\vy_j\in\overline{\Omega}_1,\\
			\infty\quad&\text{otherwise}.
		\end{cases}
	\end{equation}
\end{itemize}
Note that we have included the delay $r_1$ (which is actually zero here) as a way to streamline \cref{eq:timeE,eq:timeV} for the upcoming section. See \cref{fig:timetable} for a diagram showcasing these formulas.

The smallest entry of $\va_1$ is the time at which the first ``boundary event'' (reflection or diffraction) can happen\footnote{We say ``\emph{can} happen'' since not all vertices can cause diffraction, when hit from a certain point source. This issue is discussed in \cref{sec:diffract}.}. The index of the smallest entry tells us whether the event is a reflection (index $1\leq j\leq n_e$) or a diffraction (index $n_e+1\leq j\leq n_e+n_v$), and also what edge/vertex causes the event. From here, we use the models described in \cref{sec:reflect,sec:diffract} to build $\uapp_2$, by computing $\vxi_2$, $r_2$, $\Omega_2$, and $\zeta_2$.

Then, the second timetable vector $\va_2$ can be computed by replacing all subscripts ``1'' by ``2'' in \cref{eq:timeE,eq:timeV}. This is followed by the construction of $\uapp_3$, and so on. The process continues until all not-yet-explored entries of the timetable are larger than $T+R$, where, as mentioned above, $R$ is the half-width of the support of forcing term and initial conditions. Indeed, starting from this time instant, the would-be next terms of $\uapp$ do not affect the approximation anymore, since, due to the finite speed of wave propagation, they only act (on $\Omega$) after the end of the time horizon, i.e., for $t>T$. The total number of field components $N$ is simply the number of vectors in the timetable.

We summarize the overall procedure for the construction of the terms $\uapp_n$ in \cref{algo:timetable}. For ease of presentation, once an entry of the timetable has been explored, it is set to $\infty$ as a way for the algorithm to ignore it from that point forward.

\begin{algorithm}
	\caption{Step-by-step construction of the surrogate model}
	\label{algo:timetable}
	\begin{algorithmic}
		\STATE{Set $N\gets 1$, find $\Omega_1$ as in \cref{eq:cone}, and define $\uapp_1$ as in \cref{eq:first}}
		\STATE{Define $\va_1\in\bR^{n_e+n_v}$ using \cref{eq:timeE,eq:timeV}}
		\STATE{Set $i\gets 1$ and $j\gets\argmin_{j=1,\ldots,n_e+n_v}\left(\va_1\right)_j$}
		\WHILE{$\left(\va_i\right)_j\leq T+R$}
		\STATE{Set $\left(\va_i\right)_j\gets\infty$ and $N\gets N+1$}
		\IF{$j\leq n_e$}
		\STATE{Find $\vxi_N$, $r_N$, $\Omega_N$, and $\zeta_N$ as in \cref{sec:reflect} \hfill $\gets$ Reflection from edge $j$}
		\ELSE
		\STATE{Find vertex index $j'\gets j-n_e$}
		\STATE{Find $\vxi_N$, $r_N$, $\Omega_N$, and $\zeta_N$ as in \cref{sec:diffract} \hfill $\gets$ Diffraction from vertex $j'$}
		\ENDIF
		\STATE{Define $\uapp_N$ from $\vxi_N$, $r_N$, $\Omega_N$, and $\zeta_N$, as in \cref{eq:sum}}
		\STATE{Define $\va_N\in\bR^{n_e+n_v}$ using \cref{eq:timeE,eq:timeV}, with ``$N$'' replacing ``$1$'' in subscripts}
		\STATE{Set $(i,j)\gets\argmin_{i=1,\ldots,N, j=1,\ldots,n_e+n_v}\left(\va_i\right)_j$}
		\ENDWHILE
	\end{algorithmic}
\end{algorithm}

In \emph{trapping} domains, see, e.g., \cref{sec:numerical:room}, the number of terms $N$ might be rather large due to waves repeatedly ``bouncing back and forth'' between two or more edges/vertices. A large $N$, although necessary for a good approximation of all wavefronts, is undesirable since it increases the computational cost of both the construction of the surrogate $\uapp$ and its evaluation.

In \cite{lee_efficient_1988}, it is suggested to set a fixed upper bound (e.g., 1 or 2) on the maximum number of successive diffraction events that are taken into account. Here, we propose and employ an alternative method that allows for a finer control on whether a certain field component is included or not, based on its magnitude. Specifically, we remove all terms $\uapp_n$ that are smaller (in magnitude) than a certain tolerance $\textup{tol}$, uniformly over $\vx$ and $t$. This can be done as a post-processing step (thus speeding up the evaluation of $\uapp$ but not its construction) or even while building the surrogate itself. This can be achieved with a simple modification of \cref{algo:timetable}, by introducing a test on the magnitude of each soon-to-be-added wave contribution $\uapp_n$, discarding terms that are too small.

\section{Modeling reflection}\label{sec:reflect}
We now present the strategy for modeling reflection due to an edge $\gamma$ of the domain boundary $\partial\Omega$. We rely on the well-known ``\emph{geometrical optics}'' model, which describes wave propagation in terms of rays, not accounting for any diffraction \cite{mcnamara1990introduction}. We assume that we are adding a new term $\uapp_n$ to the surrogate model in \cref{eq:sum}, due to a reflection phenomenon caused by the field component $\uapp_i$, $i<n$. Specifically, we assume that a wave originating at $\vxi_i$ hits the edge $\gamma\subset\partial\Omega$, which, in particular, requires $\gamma\cap\overline{\Omega_i}\neq\emptyset$. We need to prescribe several ingredients.

\paragraph{Spatial correction $r_n$} We just transfer $r_n$ over from the incoming wave: $r_n=r_i$. Indeed, as we will see in \cref{sec:diffract}, we require the term $r_n$ only when modeling diffraction. 

\begin{figure}[tb]
	\centering
	\includegraphics{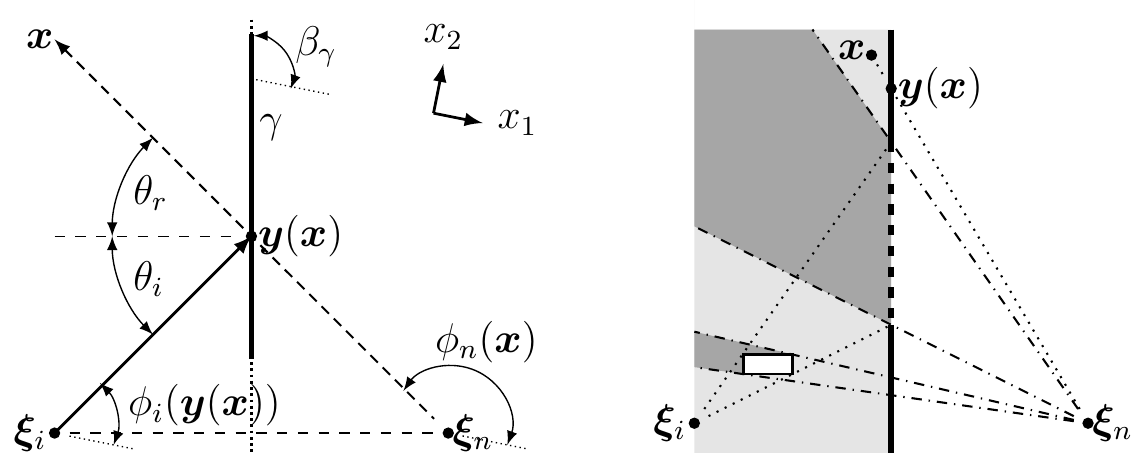}
	\caption{Graphical representation of a reflection off edge $\gamma$. On the left, the law of reflection prescribes $\theta_r=\theta_i$. We show the straight line $\widetilde{\gamma}$ supporting $\gamma$ with a dotted stroke. For a given observation point $\vx$, $\vy(\vx)$ denotes the point of incidence of the reflected field component. On the right, computation of the spatial support $\Omega_n$ (light-grey area) and its complementary shadow zone $\Omega\setminus\Omega_n$ (dark-grey area) for the reflected field component, in the presence of a rectangular obstacle. The dashed portion of edge $\gamma$ denotes the shadow $\gamma\setminus\gamma^{(i)}$. The shadow zone consists of two connected components.}
	\label{fig:reflection}
\end{figure}

\paragraph{Source point $\vxi_n$} We use the method of images, which gives the position of $\vxi_n$ as the reflection of $\vxi_i$ with respect to the edge $\gamma$:
\begin{equation}\label{eq:source_reflect}
	\vxi_n=2\argmin_{\vz\in\widetilde{\gamma}}\norm{\vz-\vxi_i}-\vxi_i,
\end{equation}
where $\widetilde{\gamma}\subset\bR^2$ is the straight line on which edge $\gamma$ lies. See \cref{fig:reflection} (left).

\paragraph{Weight function $\zeta_n$} Let $\vx-\vxi_n$ be a generic point where we wish to evaluate the weight function $\zeta_n$. We define the \emph{incidence} point $\vy(\vx)$ as the intersection (if any) between edge $\gamma$ and the segment from $\vxi_n$ to $\vx$. See \cref{fig:reflection} (left). According to the method of images, the amplitude of the reflected wave is equal (up to sign) to the amplitude of the incoming wave:
\begin{equation}\label{eq:reflectionimages}
    \zeta_n(\vx-\vxi_n)=\pm\zeta_i(\vy(\vx)-\vxi_i).
\end{equation}
Above and throughout this section, the sign $\pm$ depends on the type of boundary conditions on the edge $\gamma$ (``$+$'' for Neumann, ``$-$'' for Dirichlet).

Now, recall that we are assuming all weight functions to be positive-homogeneous.
According to \cref{rem:homogeneous}, this means that $\zeta_i(\vx-\vxi_i)$ is only a function of the angle $\phi_i(\vx)$ between $\vx-\vxi_i$ and the positive $x_1$-axis. See \cref{fig:reflection} (left) for a graphical depiction. Specifically, with an abuse of notation, let $\zeta_i(\vx-\vxi_i)=\zeta_i(\phi_i(\vx))$ and $\zeta_n(\vx-\vxi_n)=\zeta_n(\phi_n(\vx))$, where the ``new'' angle-dependent functions $\zeta_i$ and $\zeta_n$ are $2\pi$-periodic. By \cref{eq:reflectionimages}, we deduce the property
\begin{align}
\zeta_n(\phi_n(\vx))=\pm\zeta_i(\phi_i(\vy(\vx)))
=&\pm\zeta_i(2\beta_\gamma-\phi_n(\vy(\vx)))\nonumber\\
=&\pm\zeta_i(2\beta_\gamma-\phi_n(\vx)),
\end{align}
where $\beta_\gamma$ is the angle between edge $\gamma$ and the positive $x_1$-axis. This uniquely identifies $\zeta_n$, given $\zeta_i$ and $\beta_\gamma$.

\paragraph{Spatial support $\Omega_n$} We first identify what portion of $\gamma$ is actually ``lit'' by $\uapp_i$: $\gamma^{(i)}=\gamma\cap\overline{\Omega_i}$. Note that we may have $\gamma\neq\gamma^{(i)}$, for instance when obstacles are present between $\vxi_i$ and $\gamma$. See \cref{fig:reflection} (right) for an illustration. Then, roughly speaking, we define the new support $\Omega_n$ as the union of all line segments from $\vxi_n$ that pass through $\gamma^{(i)}$. To be more precise, given $\vx\in\Omega$, let $\vy(\vx)$ be the intersection (if any) between $\gamma$ and the line segment from $\vxi_n$ to $\vx$. Also, if $\vy(\vx)$ exists, we define $\tau_0(\vx)=\norm{\vy(\vx)-\vxi_n}/\norm{\vx-\vxi_n}\in(0,1)$, which satisfies $\vy(\vx)=\vxi_n+\tau_0(\vx)(\vx-\vxi_n)$. The new support is defined as
\begin{equation}\label{eq:reflcone}
	\Omega_n=\left\{\vx\in\Omega\ :\ \vy(\vx)\in\gamma^{(i)}\text{ and }\vxi_n+\tau(\vx-\vxi_n)\in\Omega\ \ \forall\tau_0(\vx)<\tau\leq 1\right\}.
\end{equation}
\begin{figure}
	\centering
	\includegraphics{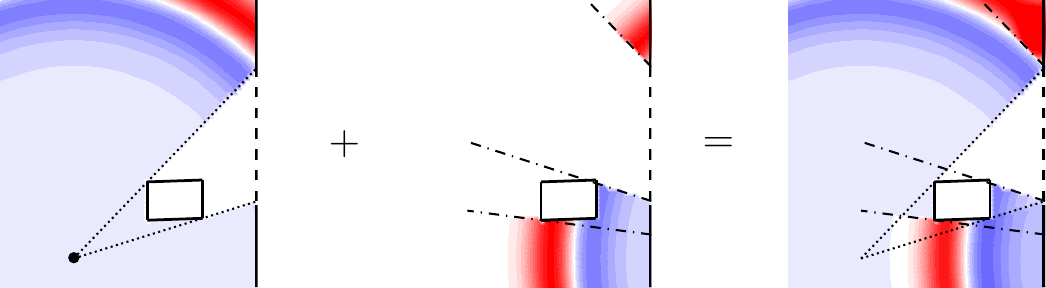}
	\caption{Example of reflection off an edge in the presence of an obstacle, from \cref{fig:reflection}. Neumann conditions are imposed on all edges. Source wave (left), reflected wave (middle), and superimposition of the two (right). Note how the obstacle creates a shadow zone for source and reflected waves. For simplicity, in this plot we are not showing any reflection or diffraction effects due to the rectangular obstacle, since they would be modeled at different stages of the algorithm.}
	\label{fig:rifl}
\end{figure}

\Cref{fig:rifl} represents a possible output of the numerical algorithm. In this case, we simulate only the reflections, thus discarding, for the time being, any effect due to diffraction. It is clear that, by modeling reflection effects only, we may obtain a discontinuous approximation of the solution of our target problem, with discontinuities arising at the boundaries of the spatial supports identified so far. As we will see in the next section, introducing diffraction in our approximation will allow us to obtain a continuous approximation $\uapp$.

\section{Modeling diffraction}\label{sec:diffract}

Here, we describe a strategy for modeling waves diffracted by a vertex of the domain boundary $\partial\Omega$. This is required in building a new field component $\uapp_n$ whenever the smallest unexplored entry of the timetable is related to a vertex, i.e., $j>n_e$ in \cref{algo:timetable}, cf.~\cref{sec:approx}.

In our modeling, we take inspiration from the uniform theory of diffraction (UTD) \cite{kouyoumjian1974uniform,mcnamara1990introduction,nethercote_analytical_2020}. UTD provides an effective description of diffraction in the frequency domain, which we need to specialize for our time-domain modeling. This task is difficult for several reasons. On the one hand, standard approaches for frequency-to-time-domain conversion involve the inversion of the Fourier transform, which is rather costly in a numerical setting, where integration/convolution must be replaced by quadrature. (This issue actually affects also time-domain diffraction modeling based on impulse responses, see, e.g., \cite{calamia_svensson_06,svensson_99}.) On the other hand, UTD is an \emph{asymptotic} method, whose accuracy relies on assuming the frequency to be sufficiently large. As such, an inverse Fourier transform is not guaranteed to provide sound results, whenever the signal bandwidth includes low frequencies.

Both these reasons behoove us to develop a novel time-domain version of UTD. Specifically, we wish our diffraction modeling to fit in the framework of \cref{eq:sum}, so that:
\begin{itemize}
    \item Diffraction is modeled as a wave outgoing from a point source at some $\vxi_n$. This is actually in line with standard approaches, e.g., \cite{biot_tolstoy,svensson_savioja_overview}, where diffraction is modeled through virtual sources on the vertex (in 2D) or edge (in 3D) that causes the diffraction event. This motivates the choice of the center $\vxi_n:=\vy_{j'}=\vy_{j-n_e}$, the diffraction vertex (we are employing the notation of \cref{algo:timetable}).
    \item Diffraction has a spatial support $\Omega_n$, which (again in accordance to diffraction theory) we define as the set of all points that are visible (along straight-line paths) from $\vxi_n$, i.e.,
    \begin{equation}\label{eq:diffcone}
	\Omega_n=\left\{\vx\in\Omega\ :\ \vxi_n+\tau(\vx-\vxi_n)\in\Omega\ \ \forall0<\tau\leq 1\right\}.
    \end{equation}
    \item The space-time dependence of the diffraction amplitude should be separable into an angular component $\zeta_n$, and a radial-temporal component, which, in fact, is assumed to coincide with the free-space wave propagation profile $\Psi$.
\end{itemize}

The last step is the most critical one, as it involves a simplifying modeling choice, specifically made to fit the \emph{approximate} diffraction wave within our framework. Still, there are strong theoretical foundations for this, as we proceed to explain.

\subsection{Modeling the angular modulation}\label{sec:diffract:mod}
From a phenomenological point of view, the angular modulation $\zeta_n$ can be related to the \emph{diffraction coefficient} appearing in geometrical diffraction theory \cite{keller1962geometrical,mcnamara1990introduction}. Indeed, the diffraction coefficient $D$ has exactly the desired role of scaling factor for the diffraction amplitude, which depends on the angle that the target measurement point $\vx$ forms with the edges adjacent the diffraction vertex. As such, we set $\zeta_n$ equal to the UTD diffraction coefficient $D$.

We define $D$ in full detail in \cref{app:diffract}. Here, for our discussion, we only need the following basic facts: in 2-dimensional UTD, the diffraction coefficient $D$ depends on the above-mentioned angle between $\vx$ and $\vxi_n$, but also on (i) the domain geometry locally around the diffraction point, (ii) the location of the source that causes the diffraction, (iii) the frequency of the incident signal, and (iv) the distance $\norm{\vx-\vxi_n}$.

The last two above-mentioned items are problematic in our framework: since we work in the time domain, we do not have a single incident frequency $k$; moreover, due to our separability assumption, we require the angular modulation to be independent of $\norm{\vx-\vxi_n}$. However, frequency and radial distance always appear together, as a product, in the UTD diffraction coefficient. Specifically, such product $\ks:=k\norm{\vx-\vxi_n}$ can be interpreted as the distance between $\vx$ and the diffraction point $\vxi_n$, measured in wavelength units.

In practice, the value of $\ks$ determines how quickly the magnitude of $D$ decays to $0$, as one moves away from a so-called ``shadow boundary'', i.e., the boundary of the spatial support of either the incident or the reflected wave(s). Equivalently, $\ks$ determines the width of the \emph{transition region}, which encompasses any shadow boundaries. See \cref{app:diffract,fig:diff_coeff}. Also, note that UTD requires $\ks$ to be large enough ($\geq 1$), and that $\ks\to\infty$ yields the geometrical-optics setting, i.e., a diffraction-free model, with $\zeta_n\equiv 0$.

Since the dependence of $D$ on the two troublesome terms happens only through $\ks$, we propose the following strategy for defining a $k$- and $\norm{\vx-\vxi_n}$-independent angular modulation $\zeta_n$: we set $\zeta_n$ as the diffraction coefficient $D$ obtained for some value of $\ks=\overline{\ks}$ that is fixed \emph{a priori}. Such value $\overline{\ks}$ should be sufficiently large not to cause issues in UTD, but also sufficiently small to avoid the geometrical-optics pitfall.

Obviously, the specific choice of $\overline{\ks}$ should be based on the ultimate approximation target: minimizing the discrepancy between $u$ and $\uapp$. However, it turns out that $D$ depends relatively mildly on the value of $\ks$. See \cref{app:diffract} for a theoretical justification, and \cref{sec:numerical:wedge} for some numerical tests. In our experiments, we settle for the value $\overline{\ks}=10$. A more careful quantitative investigation of the role of $\ks$ on the diffraction approximation accuracy is envisioned as a future research direction.

Instead of using frequency-domain UTD, the time-domain diffraction coefficient could be computed by convolution of the incoming wave with a ``diffraction impulse response'', as done in \cite{svensson_99,torres_computation_2001}. Such a convolution would need to be approximated (by quadrature) over and over, at each diffraction event. Notably, since the diffraction impulse response is \emph{not} a Dirac delta, a convolution with it naturally introduces a radial modulation too. This entails that the diffraction wave does not behave like an angular modulation of the free-space wave $\Psi$, as assumed in our model.

Our modeling of diffraction through a diffraction coefficient that does not depend on time nor radius, although biased in the above-described way, is very efficient, as it avoids the expensive computation of convolutions with the diffraction impulse response. As we will show in our numerical experiments, this increased efficiency comes at the cost of a relatively small modeling error. As such, our approach remains justified.

\subsection{Estimating the modeling error}
We can assess the quality of our diffraction model \emph{a priori}, by comparing its resulting approximation of the diffraction wave ($\uapp_n$) with that obtained via the above-mentioned convolution approach \cite{svensson_99}, which we denote by $\mathring{u}_n$ here.

As before, assume that field component $\uapp_i$ impinges on vertex $\vxi_n$, causing a diffraction event emanating from there. Given an arbitrary space-time position $(\vx,t)$, the local discrepancy between the diffraction waves obtained with the two methods can be expressed by inverse Fourier transform, as
\begin{align*}
    \uapp_n(\vx,t)-\mathring{u}_n(\vx,t)=&\widetilde{D}_{\overline{\mu}}\int_{\mathbb{R}}\mathcal{F}[\uapp_i](\vxi_n,k)e^{ikt}\textrm{d}k-\int_{\mathbb{R}}D_k\mathcal{F}[\uapp_i](\vxi_n,k)e^{ikt}\textrm{d}k\\
    =&\int_{\mathbb{R}}\left(\widetilde{D}_{\overline{\mu}}-D_k\right)\mathcal{F}[\uapp_i](\vxi_n,k)e^{ikt}\textrm{d}k.
\end{align*}
Above, $\mathcal{F}$ denotes the Fourier transform operator, whereas $D_k=D_k(\vx-\vxi_n)$ is the UTD diffraction coefficient, which depends on both the position $\vx$ and the wavenumber $k$, cf.~\cref{app:diffract}. On the other hand, $\widetilde{D}_{\overline{\mu}}$ is the diffraction coefficient used in our modeling, which is $k$-independent and depends on $\vx$ only through the direction of $\vx-\vxi_n$. (Note that, in geometrical optics, one would set $\widetilde{D}_{\overline{\mu}}=0$.)

It follows that, to achieve a small error, one should pick a value of $\widetilde{D}_{\overline{\mu}}$ close to that of $D_k$, specifically at wavenumbers $k$ that lie closest to the Fourier spectrum of $\uapp_i(\vxi_n,\cdot)$. In this way, we can expect that our approach, despite being based on a heuristic choice of $\overline{\mu}$, independently of $k$, should do significantly better than geometrical optics.

This being said, our method has limitations, which are better seen by looking at the following global-in-wavenumber bound on the error magnitude:
\begin{equation}\label{eq:diffr_error_bound}
    \abs{\uapp_n(\vx,t)-\mathring{u}_n(\vx,t)}\leq\sup_k\abs{\widetilde{D}_{\overline{\mu}}-D_k}\norm{\mathcal{F}[\uapp_i](\vxi_n,\cdot)}_{L^1(\mathbb{R})}.
\end{equation}
Our method may lead to gross errors if $D_k$ has large variations in $k$. In this respect, it is worth nothing that, thanks to the availability of explicit formulas for $D_k$, cf.~\cref{app:diffract}, it is possible to obtain an upper-bound for $|D_k-\widetilde{D}_{\overline{\mu}}|$ by direct inspection. Still, bound \cref{eq:diffr_error_bound} is pessimistic, as we will show in our numerical experiments.

As a practical alternative, we propose a heuristic \emph{a posteriori} error indicator, with the objective of obtaining better (but still approximate) information on the magnitude of the diffraction modeling error. To this aim, we apply the following multi-fidelity argument. Let $\uapp$ be the approximation obtained through our method, including also diffraction terms, modeled as in \cref{sec:diffract:mod}. On the other hand, let $\uapp_{\text{GO}}$ be the approximation obtained with our approach, \emph{disregarding diffraction completely} (``GO'' stands for ``geometrical optics''). As indicator of the approximation error $\abs{\uapp-u}$, we simply use the cross-fidelity discrepancy $|\uapp_{\text{GO}}-\uapp|$. The justification behind this follows.

First, we note that, since reflection effects are modeled exactly in our approach, cf.~\cref{sec:reflect}, $u-\uapp_{\text{GO}}$ coincides with the superposition of all diffraction-caused effects. Second, we make the heuristic assumption that diffraction effects in $\uapp$ are approximated with at most a $50\%$ relative error. (Of course, we have no way of verifying this in practice, as we would need access to the exact solution.) Accordingly, we have that $2|\uapp-u|\leq|\uapp_{\text{GO}}-u|$ and, by the triangular inequality,
\begin{equation*}
    \abs{\uapp-u}\leq\abs{\uapp_{\text{GO}}-u}-\abs{\uapp-u}\leq\abs{\uapp_{\text{GO}}-\uapp}.
\end{equation*}
We assess the reliability of this estimate in our numerical experiments in the next section.

\section{Numerical results}\label{sec:numerical}
In our experiments, we require a ``reference'' solution of \cref{eq:wave} to validate our results. To this effect, we use the solution $\ufem$ obtained by discretizing \cref{eq:wave} with:
\begin{itemize}
	\item the P1-finite element (FE) method with mass-lumping, over a regular triangulation (mesh) of the physical domain $\Omega$;
	\item explicit leapfrog timestepping with a uniform time step that satisfies the CFL condition on the chosen mesh.
\end{itemize}
See \cite{Cohen2008,HairerLubich} for more details on this discretization strategy.

If the domain $\Omega$ is unbounded, we first need to truncate it in such a way that reflections from the non-physical truncation boundary do not affect the solution in the region of interest for $t<T$. Recalling that the problem data are supported in a ball of radius $R$ and center $\vzero$, this can be done, e.g., by truncating $\Omega$ at the sphere with radius $R+T$ (we recall that we are assuming a unit wave speed) and center $\vzero$. In our tests, we rely on FEniCS \cite{fenics} to carry out the FE discretization on 2-dimensional domains $\Omega$. Note that, instead, unbounded geometries are allowed in our proposed approach, making domain truncations unnecessary.

All our tests are performed in Python 3.8 on a machine with an 8-core 4.70 GHz Intel\textsuperscript{\textregistered} processor. For reproducibility, our code is made available at \url{https://github.com/pradovera/ray-wave-2d}.

\subsection{Some simple wedges}\label{sec:numerical:wedge}
As a way to assess our proposed method in simple settings, we consider four different ``wedge'' domains. We define $\Omega$ to be one of the portions of the plane $\bR^2$ delimited by two straight lines intersecting at a point, with $\alpha$ being the outer angle at such point. The specific choices of wedge angles $\alpha$ are reported in \cref{tab:wedge} for the four cases.

\begin{figure}[tb]
	\centering
	\includegraphics{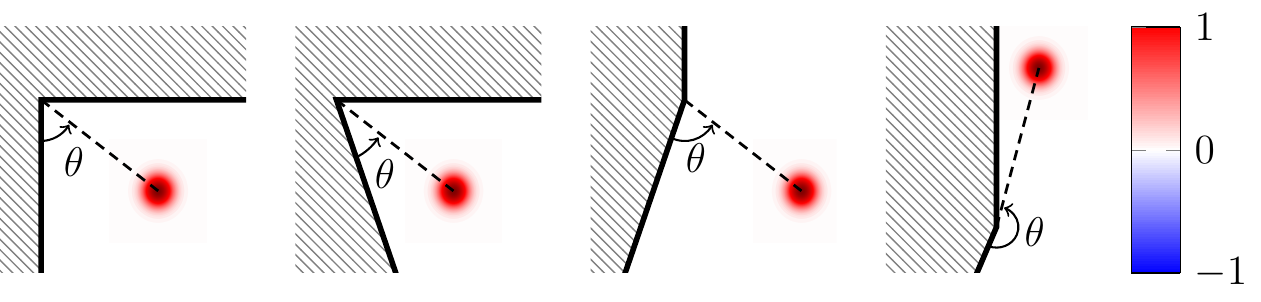}
	\caption{Initial conditions for the wedge examples, indexed \#1 through \#4 from left to right. The (dashed) distance between the center of the Gaussian and the boundary vertex is 4 units in all cases.}
	\label{fig:iwedge}
\end{figure}

\begin{table}[tb]
	\centering
	\begin{scriptsize}
	\begin{tabular}{c||c|c||c|c}
		Example & exterior & incidence & number $N$ of & \multirow{2}{*}{$\sup_{t\in[0,T]}\textup{err}(t)$, cf.~\cref{eq:err}}\\
		index & angle $\alpha$ & angle $\theta$ & components $\uapp_n$ & \\
		\hline
		\#1 & $1.5\pi$ & $0.313\pi$ & 5 & $0.146\%$ \\
		\#2 & $1.62\pi$ & $0.192\pi$ & 8 & $0.828\%$ \\
		\#3 & $0.879\pi$ & $0.434\pi$ & 4 & $3.14\%$ \\
		\#4 & $0.879\pi$ & $1.04\pi$ & 3 & $4.68\%$ \\
	\end{tabular}
	\end{scriptsize}
	\caption{Setup for the four wedge examples. The angle $\theta$ is as in \cref{fig:iwedge}. The last two columns refer to the results of our algorithm.}
	\label{tab:wedge}
\end{table}

We set up a wave-propagation problem like \cref{eq:wave}, with $u_0$ an isotropic Gaussian with standard deviation $0.2$. The center of $u_0$ is at a point located at a 4-unit distance from the wedge vertex, in the direction determined by the ``incidence angle'' $\theta$. See \cref{fig:iwedge} for a representation of the initial conditions in the four cases. We set $u_1=f=0$, we enforce Neumann boundary conditions on the whole $\partial\Omega$, and we seek the solution at the final time $T=5$, i.e., $1$ time unit after the wave crest has reached the wedge vertex.

\begin{figure}[tb]
	\centering
	\includegraphics{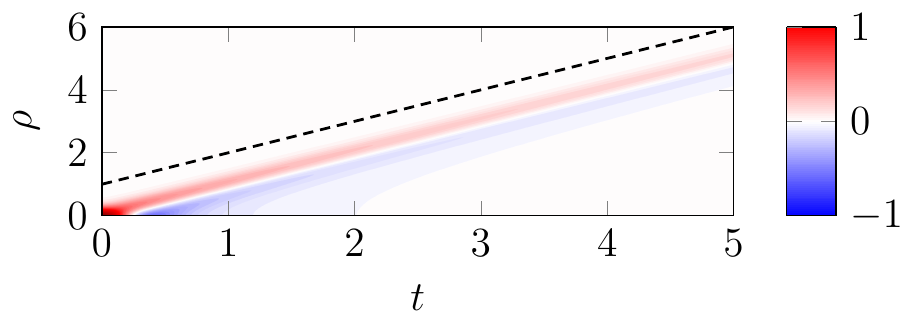}
	\caption{Free-space solution $\Psi$. The dashed line denotes the upper bound of the ``causality cone'' of $\Psi$, i.e., $\rho=t+R$, with $R=1$.}
	\label{fig:freespace}
\end{figure}

\begin{figure}[tb]
	\centering
	\includegraphics{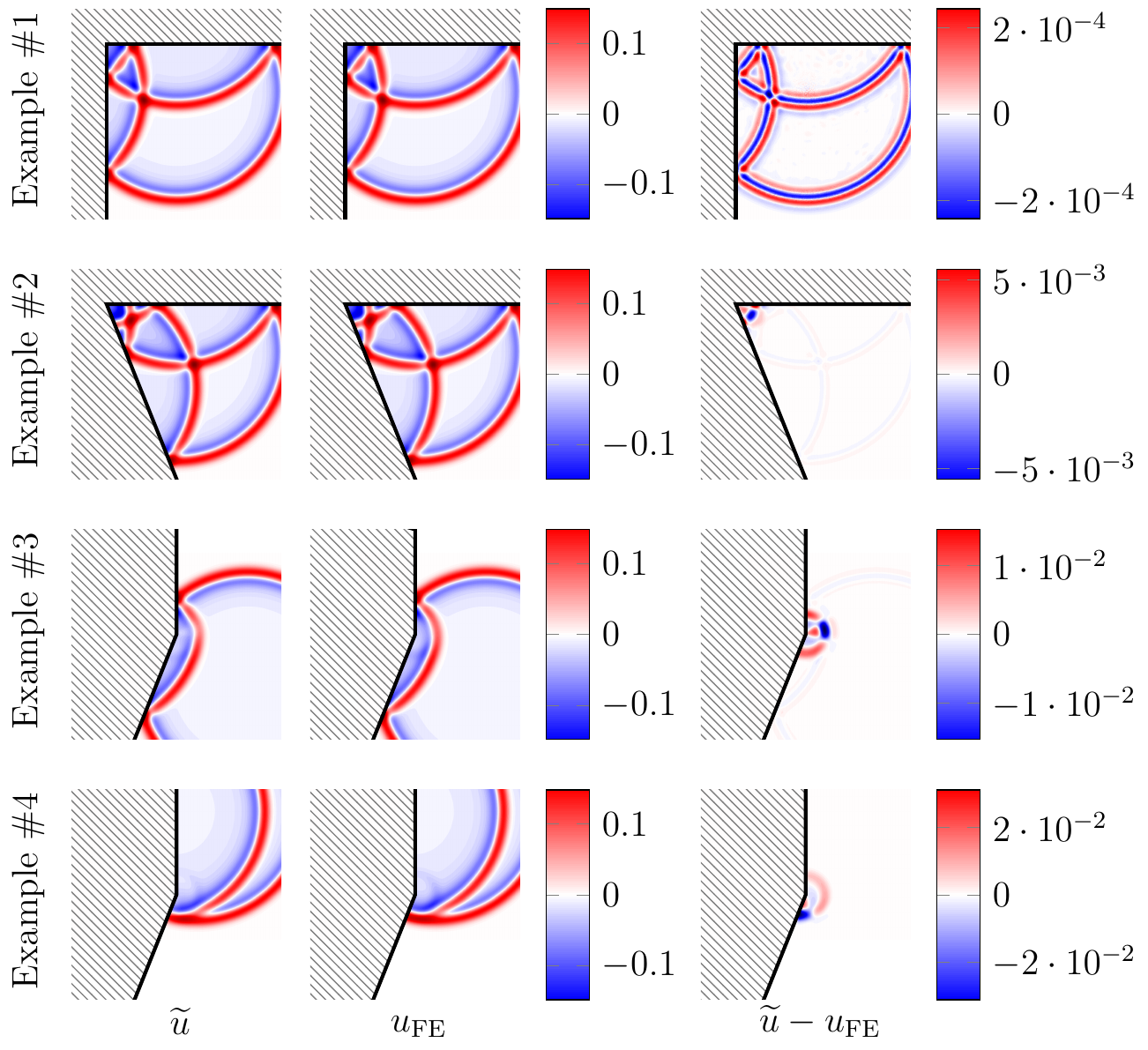}
	\caption{Results for the four wedge examples. Each row pertains to a different example. In each row, from left to right: surrogate solution, FE solution, and error. The color scales for the first two columns are the same. All results are shown at the final time $t=T$.}
	\label{fig:wedge}
\end{figure}

To this aim, we employ our proposed approach, see \cref{sec:approx}. First, we compute an approximation of the free-space solution $\Psi$, which solves \cref{eq:Psi}, by employing the P1-FE method with explicit leapfrog timestepping. Note that, since \cref{eq:Psi} is cast in polar coordinates, we only need to discretize a 1D interval with P1-FE. Since the initial condition $u_0$ is supported within the unit disk, we have $R=1$, and it suffices to approximate $\Psi(\rho,t)$ for $(\rho,t)\in[0,T+R]\times[0,T]$. Since this space-time domain is only 2-dimensional, we can afford a very fine discretization. In our experiments, we employ a $1001\times 2001$ uniform Cartesian space-time grid, i.e., the mesh size is $\delta x=\frac{T+R}{1000}$ and the time step is $\delta t=\frac{T}{2000}$. This satisfies the CFL condition. We show the resulting $\Psi$ (which, in fact, we should denote by $\Psi_{\textup{FE}}$) in \cref{fig:freespace}. Note that, since we are solving the wave equation in 2D, the magnitude of $\Psi$ decays like $\mathcal{O}(t^{-1/2})$ as $t\to\infty$, which corresponds to a spatial decay of $\mathcal{O}(\rho^{-1/2})$. In 3D, the decay would be quicker, namely, $\mathcal{O}(t^{-1})$ in time and $\mathcal{O}(\rho^{-1})$ in space.

After this preliminary step, we use the timetable-based strategy from \cref{sec:approx} to identify reflection and scattering effects, which are then added up to give the final approximation $\uapp$. We show the resulting $\uapp(\cdot,T)$ in \cref{fig:wedge}. In this figure, we also display a reference solution $\ufem(\cdot,T)$, which we obtain by direct discretization of \cref{eq:wave} by P1-FE and leapfrog timestepping, as described at the beginning of \cref{sec:numerical}.

In all four examples, we see that $\uapp$ and $\ufem$ seem qualitatively close. Notably, we can observe a good representation of the most prominent wavefronts, which are due to propagation of either the main ``free-space'' wave or to its reflections. Indeed, those wave contributions are reconstructed exactly: the only errors are the ones due to FE approximation and timestepping, which affect both $\ufem$ and $\uapp$ (the latter through the approximation of $\Psi$). Instead, some differences are present when comparing diffraction effects, which arise as circular waves about the wedge vertex. We can quantitatively observe this in the last column of both \cref{tab:wedge,fig:wedge}.

In example \#1, we observe a very small error, which, in fact, is simply the (FE and timestepping) discretization error. This is related to the fact that the \emph{wedge index} $\nu=\pi/(2\pi-\alpha)=2$ is an integer, thus making diffraction unnecessary in approximating the wave $u$ \cite{kouyoumjian1974uniform}.

In the other examples, diffraction effects are necessary to correctly identify $u$. While a good qualitative behavior can be observed in \cref{fig:wedge}, we can see in \cref{tab:wedge} that a modest error is present. Specifically, as error measure, we use the supremum over $t\in[0,T]$ of the relative $L^2(\Omega)$-approximation error, defined as
\begin{equation}\label{eq:err}
	\textup{err}(t):=\left(\int_\Omega\left(\uapp(\vx,t)-\ufem(\vx,t)\right)^2\textup{d}\vx~\bigg/\int_\Omega\ufem(\vx,t)^2\textup{d}\vx\right)^{1/2}.
\end{equation}

We see the largest error in example \#4, where the relative $L^2(\Omega)$-ap\-pro\-xi\-ma\-tion error amounts to about $5\%$. This example corresponds to a case of ``almost grazing'' incidence, with the source point being located rather close to one of the wedge's edges. In turn, this leads to two shadow boundaries located very close to each other, with overlapping transition regions, cf.~\cref{app:diffract}. We interpret this result as evidence of the fact that our diffraction model from \cref{sec:diffract} is less accurate with grazing than with non-grazing incidence. Despite this, the error remains small in all the above tests, and our diffraction model may be considered satisfactory overall.

\subsubsection{Building a cavity out of wedges}
\begin{figure}[tb]
	\centering
	\includegraphics{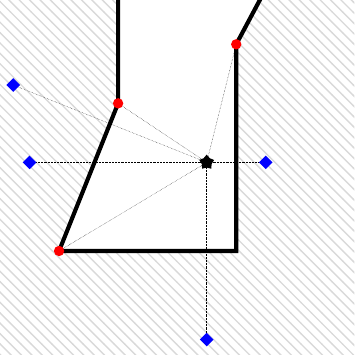}\hspace{1mm}%
	\includegraphics{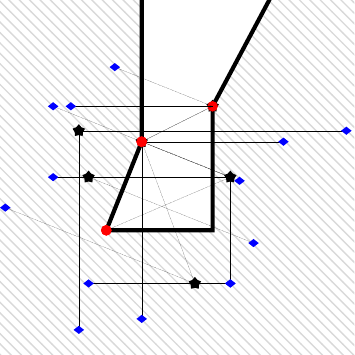}\hspace{1mm}%
	\includegraphics{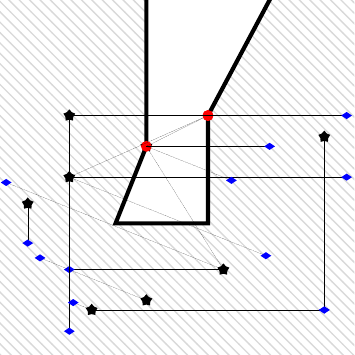}
	\caption{Graph view of the field components for the cavity example. Markers are used to denote source points $\vxi_n$. Dotted line segments denote causality effects: a black star is linked to a blue square (resp., red circle) if the field component emanating from the former causes a reflection (resp., diffraction) emanating from the latter. To avoid clutter, the field components have been separated into three plots according to causality: in the left plot we show all $7$ events (reflections or diffractions) directly caused by $\uapp_1$, in the middle plot we show all $20$ events directly caused by the $7$ field components in the left plot, etc.}
	\label{fig:tree}
\end{figure}

As a slightly more complicated example, we now combine the four wedges from the previous section to obtain the open cavity represented in \cref{fig:tree}. In this case, more reflection and diffraction effects arise, due to the trapping nature of the domain. Our initial conditions and forcing term are the same as before, but now all edges are sound-soft. Accordingly, we model them using Dirichlet boundary conditions. The time horizon is $T=9$.

Using our strategy from \cref{sec:approx}, we build the approximation $\uapp$, which contains $45$ wave terms ($1$ source wave, $32$ reflected waves, and $12$ diffraction waves). As visual reference, we display the source points of all field components $\uapp_n$ in \cref{fig:tree}. Each dotted line segment is a ``causality link'', showing how one field component causes the next, by either reflection or diffraction.

\begin{figure}[tb]
	\centering
	\includegraphics{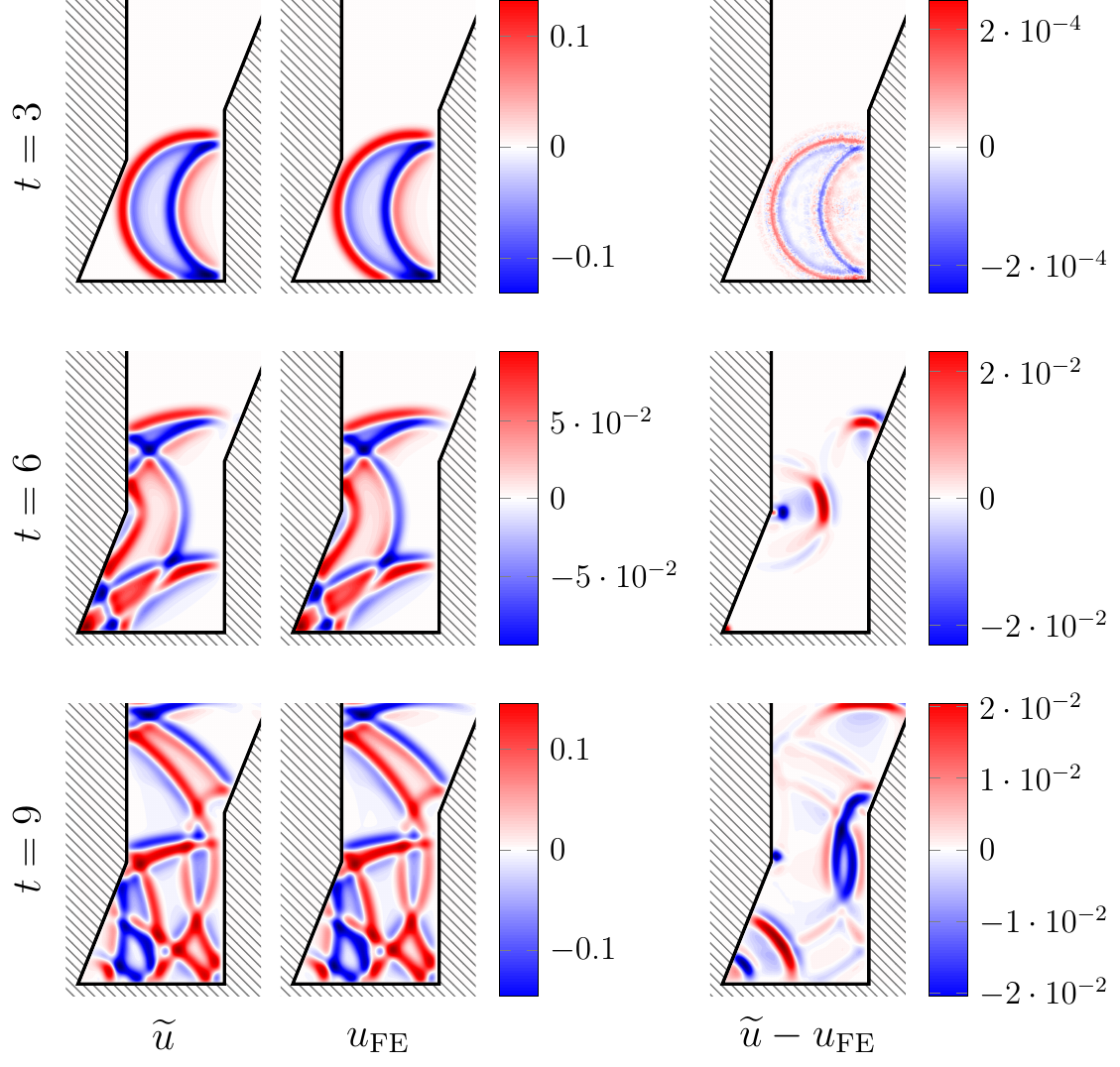}
	\caption{Results for the cavity domain. Each row corresponds to a different time instant $t\in\{3,6,9\}$, from top to bottom. In each row, from left to right: surrogate solution $\uapp(\cdot,t)$, FE solution $\ufem(\cdot,t)$, and error $\uapp(\cdot,t)-\ufem(\cdot,t)$. The color scales for the first two columns are the same.}
	\label{fig:cavity}
\end{figure}

In \cref{fig:cavity}, we compare the approximation $\uapp$ with the reference solution $\ufem$, obtained as described at the beginning of \cref{sec:numerical}, on a spatial mesh with around $10^6$ degrees of freedom. Once more, we see a good qualitative agreement between $\uapp$ and $\ufem$, with the most important features of $u$ being identified well. We see ``error waves'' of small amplitude propagating from the 3 vertices of $\Omega$ that generate diffraction effects. These correspond to errors in diffraction modeling. More quantitatively, the relative $L^2(\Omega)$-approximation errors are roughly $0.05\%$ (at $t=3$), $7.21\%$ (at $t=6$), and $12\%$ (at $t=9$).

The approximation $\uapp$ is built in approximately $500$ms, out of which $420$ms are used to compute the free-space solution $\Psi$. This is a surprisingly small time, compared to the computation of the reference $\ufem$, which takes about a minute, with the 2D meshing alone taking around $20$s. Note that we are including also the meshing time in order to have a fair comparison of the costs incurred by the two considered methods, with the target being obtaining a reasonably accurate approximate solution of the wave equation.

As an additional experiment to validate our diffraction modeling, we study how the approximation error depends on the choice of $\overline{\ks}$, which affects the diffraction coefficient by being (roughly) inversely proportional to the width of the transition regions around shadow boundaries, cf.~\cref{sec:diffract}. We compute the $L^2(\Omega)$-approximation error at times $t=6$ and $t=9$, for values of $\overline{\ks}$ between $1$ and $100$. Note that every value of $\overline{\ks}$ involves the computation of a new diffraction coefficient $D$ for each diffraction event, but otherwise does not require retraining the surrogate model $\uapp$: only the angular weights $\zeta_n$ vary, while all other terms in \cref{eq:sum} remain the same.

\begin{figure}[tb]
	\centering
	\includegraphics{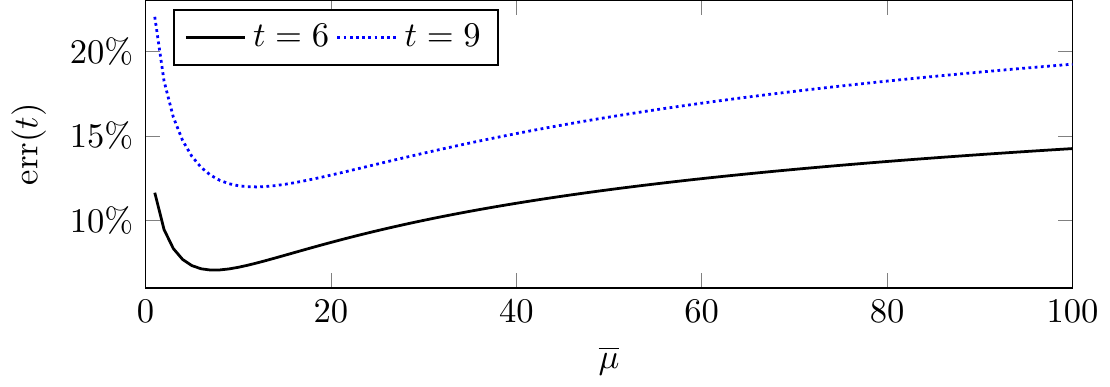}
	\caption{Relative $L^2(\Omega)$-approximation error at $t=6$ and $t=9$ for different values of $\overline{\ks}$.}
	\label{fig:ks}
\end{figure}

In \cref{fig:ks}, we show how the errors have a minimum around $\overline{\ks}=7$ for $t=6$, and around $\overline{\ks}=11$ for $t=9$. This shift in the minimum point seems to suggest that larger values of $\overline{\ks}$ should be used for approximations at longer times, an observation that may be especially relevant for simulations over long time horizons. That being said, the error appears to vary only mildly around the locations of such minima: for instance, at $t=6$, the relative error obtained for $\overline{\ks}=10$ is $7.21\%$, only slightly more than the minimum error, $7.06\%$. (Note also the narrowness of the vertical scale of the plot.) This provides an empirical justification for choosing a single value of $\overline{\ks}$ (around $10$) for all diffraction events, and for all times $t$.

\subsection{A tall room}\label{sec:numerical:room}
We now move to a simplified sound propagation problem in a room. For simplicity, we consider a 2-dimensional problem, thus assuming an infinitely tall room, and modeling line sources in the $z$-direction (e.g., an array of loudspeakers) as point sources.

\begin{figure}[tb]
	\centering
	\includegraphics{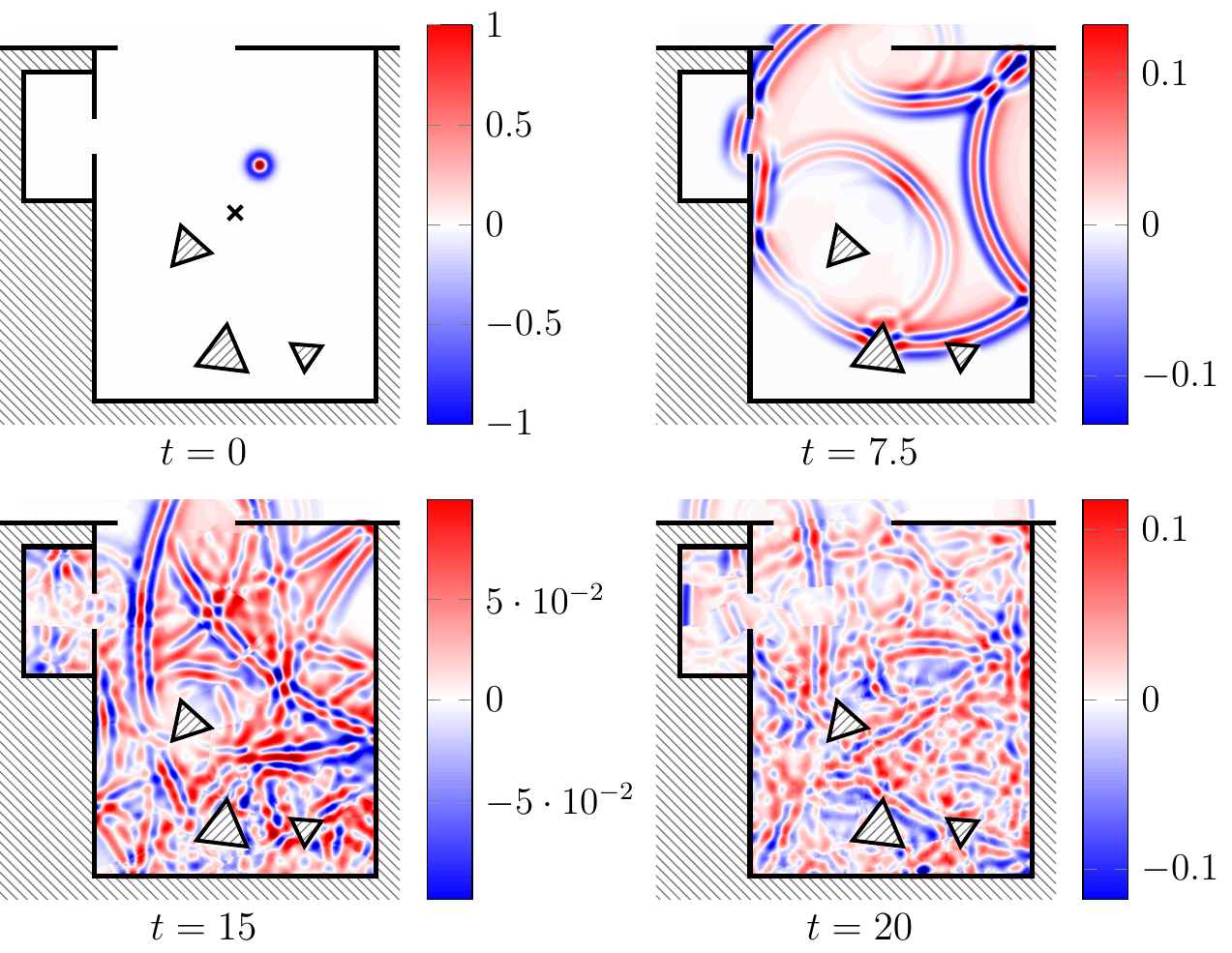}
	\caption{2-dimensional domain $\Omega$ modeling a room. Top left plot: initial condition $\uapp(\cdot,0)=u(\cdot,0)=u_0$, a Ricker wavelet; we also show the point $\vx_{\textup{trace}}$ as a cross. Top right plot: intermediate solution $\uapp(\cdot,7.5)$. Bottom left plot: intermediate solution $\uapp(\cdot,15)$. Bottom right plot: final solution $\uapp(\cdot,20)$.}
	\label{fig:room:domain}
\end{figure}

The complicated domain $\Omega\subset\bR^2$ is depicted in \cref{fig:room:domain}. It is composed of two communicating ``rooms'' with sound-hard walls, as well as of a third large room (above), which is modeled as infinitely large. In the main room, three sound-soft triangular obstacles are also present.

Setting once more $u_1=f=0$, we are interested in modeling the propagation of an initial condition $u_0$ modeled as a Ricker wavelet centered at $\vzero$, see \cref{fig:room:domain} (top left), over the time horizon $t\in[0,T]$, with $T=20$. To this aim, we employ our proposed method from \cref{sec:approx}.

As in the previous example, we start by computing an approximation of the free-space solution $\Psi=\Psi(\rho,t)$ for $(\rho,t)\in[0,T+R]\times[0,T]$, see \cref{eq:Psi}, with $R$ being the radius of the support of the initial condition $u_0$. Again, we use P1-FE with leapfrog timestepping for this.

Since many reflective surfaces face each other, the domain $\Omega$ is trapping. Accordingly, we expect the number $N$ of waves in the approximation $\uapp$ to be rather large. In the interest of reducing the number of such terms, we employ the on-the-fly tolerance-based strategy described in \cref{sec:skeleton}, removing all wave terms $\uapp_n$ whose magnitude is smaller than $\textup{tol}=2.5\cdot10^{-2}$. After this, $N=798$ terms are left. Although this value of $N$ may seem large, the evaluation of the corresponding surrogate $\uapp$ is rather quick, due to the explicit nature of each wave contribution (and to the fact that their supports are smaller than the whole $\Omega$).

\begin{figure}[tb]
	\centering
	\includegraphics{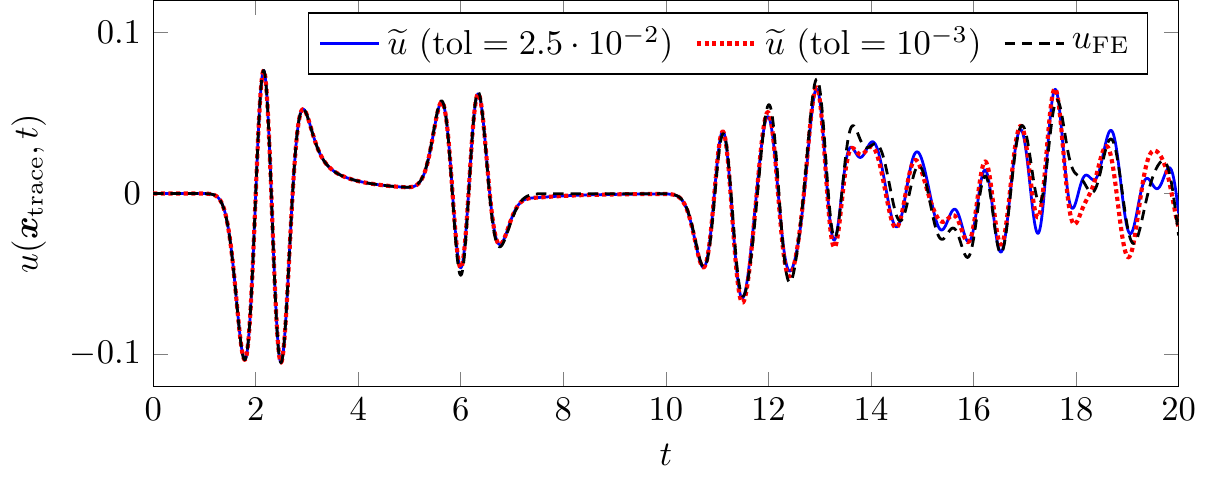}
	\caption{Value of solution at point $\vx_{\textup{trace}}=(-1,-2)$.}
	\label{fig:room:trace}
\end{figure}

We show the resulting $\uapp(\cdot,t)$ for the four times $t\in\{0,7.5,15,20\}$ in \cref{fig:room:domain}. There, we can see why so many terms are necessary for the approximation of $u$: we must model many reflection and diffraction effects. Since energy escapes the system only through the top ``door'', the wave will persist for quite a long time. Accordingly, a larger $T$ will make a larger $N$ necessary.

In order to better inspect this effect, we show the trace of the solution at the arbitrarily chosen point $\vx_{\textup{trace}}=(-1,-2)$ in \cref{fig:room:trace}. We notice that oscillations persist for $t>10$. We use this last plot also to validate our results. To this aim, we compare:
\begin{itemize}
	\item The surrogate $\uapp$ obtained as described above, with $\textup{tol}=2.5\cdot10^{-2}$.
	\item The surrogate $\uapp$ obtained with our strategy, but with $\textup{tol}=10^{-3}$. This leads to an increased number of rays $N\approx 9.6\cdot 10^3$.
	\item The reference solution $\ufem$ obtained by the P1-FE with leapfrog time\-step\-ping, as described at the beginning of \cref{sec:numerical}. The mesh size must be chosen small enough to well resolve both the initial condition and the domain $\Omega$, as well as propagation of the solution itself. In our case, we found that a mesh with approximately $1.4\cdot 10^6$ elements leads to sufficient accuracy for a comparison with both above surrogate models $\uapp$. To satisfy the CFL condition on this mesh, we choose a time step $\Delta t\approx 7\cdot 10^{-3}$.
\end{itemize}

We can observe that the two surrogates obtained with our approach give very similar results. Indeed, the cutoff tolerance $\textup{tol}$ affects the results only for sufficiently large $t$, due to the accumulation of ``small'' waves that are excluded from the coarser surrogate but included in the finer one.

Moreover, taking the FE solution as reference, we see that most of the peaks of the surrogates are aligned with the FE ones (i.e., the ``phase'' of the wave is well approximated), but there are some noticeable discrepancies in their amplitudes. This is due to the fact that, in our approach, reflection is modeled exactly, whereas the magnitudes of the diffraction waves are only approximated. For this reason, we should not expect the amplitude error to get smaller if we reduce $\textup{tol}$. The only effective way of improving the approximation would be using a more accurate diffraction model.

As a final result, we also report:
\begin{itemize}
	\item The ``construction'' time, i.e., the time required to compute the numerical solution. For $\uapp$, this means executing \cref{algo:timetable}. For $\ufem$, this means building the mesh, assembling the FE stiffness and (lumped) mass matrices, and carrying out the timestepping.
	\item The ``evaluation'' time, i.e., the time required to evaluate the numerical solution ($\uapp$ or $\ufem$) at a single $(\vx,t)$-point.
\end{itemize}
They can be found in \cref{tab:room}.

\begin{table}[tb]
	\centering
	\begin{tabular}{c|c|c|c}
		 & $\uapp$ ($\textup{tol}=2.5\cdot10^{-2}$) & $\uapp$ ($\textup{tol}=10^{-3}$) & $\ufem$ \\
		\hline
		Construction & $1.03\cdot 10^1$ & $1.39\cdot 10^2$ & $9.21\cdot 10^1$ \\
		Evaluation & $3.38\cdot 10^{-4}$ & $4.71\cdot 10^{-3}$ & $3.89\cdot 10^{-5}$ \\
	\end{tabular}
	\caption{Computational times (in seconds) for the room test case. To obtain more statistically significant results, each displayed time is the average over $3$ (resp. $10^3$) runs of the construction (resp. evaluation) phase with identical parameters.}
	\label{tab:room}
\end{table}

We can observe the increased construction and evaluation times that result from decreasing $\textup{tol}$. Moreover, we see that, in this example, our proposed approach is more competitive in the construction phase, but less so in the evaluation phase. The larger evaluation time of our surrogate is ultimately due to the nonlinearity of the functions that appear in $\uapp$. In this context, it may be surprising to note that, in our implementation, the most expensive step in evaluating $\uapp$ (taking about half of the online time) is determining whether an evaluation point is in the spatial supports $\Omega_n$ or not. On the other hand, evaluating the FE solution at a space-time point is an extremely cheap operation, essentially corresponding to a vector dot product.

This being said, we can notice a clear advantage of our approach in this context. Once the approximation $\uapp$ has been built, it can easily be evaluated at \emph{arbitrary} locations in space \emph{and time}, namely, not only at the final time $t=T$. This is because all the terms defining the approximate expansion \cref{eq:sum} ($\Psi$, $\vxi_n$, etc.) are cheap to store in memory. One cannot usually do the same with the FE solution. Indeed, once the time-stepping has been carried out, the evaluation of $\ufem$ at arbitrary instants before $t=T$ (by space-time interpolation) is possible only if the whole sequence of solution at all time-steps has been stored, a practice that is commonly avoided due to the (often) unfeasible memory requirements.

\subsubsection{A time-harmonic source}\label{sec:numerical:helm}
One of the advantages of our approach is that it allows changing the source terms of the problem in a seamless way. Notably, under minor technical constraints (e.g., the spatial support of the new source term should not be larger than the old one), this kind of change does not require training a new surrogate.

To showcase this, we approximate the wave propagating from a time-harmonic point source at $\vx=\vzero$ with angular frequency $\omega>0$. In our tests, we pick $\omega\in\{2\pi,10\pi\}$. To this aim, we define $u$ as the solution of the following ($\omega$-dependent) problem:
\begin{equation}\label{eq:helmholtz}
	\begin{cases}
		\partial_{tt}u(\vx,t)=\Delta u(\vx,t)-\omega^2\sin(\omega t)g(\vx)\quad&\text{for }(\vx,t)\in\Omega\times(0,T),\\
		u(\vx,0)=0\quad&\text{for }\vx\in\Omega,\\
		\partial_tu(\vx,0)=0\quad&\text{for }\vx\in\Omega,\\
		\partial_\nu u(\vx,t)=0\quad&\text{for }(\vx,t)\in\partial\Omega\times(0,T],
	\end{cases}
\end{equation}
where $g$ denotes a narrow 2-dimensional Gaussian, centered at $\vzero$ and with standard deviation $0.05$.

As usual, we define $\Psi=\Psi(\rho,t)$ as the ($\omega$-dependent) solution of the free-space version of \cref{eq:helmholtz} in radial-temporal coordinates. To obtain an approximation for the wave $u$ generated by the time-harmonic source for an arbitrary $\omega$, it suffices to plug the corresponding $\Psi$ in each term of the surrogate $\uapp$ from the previous section! We show the results of our approximation in \cref{
fig:room:h:trace}, where we can observe a very good agreement between approximation and reference FE solution.

\begin{figure}[tb]
	\centering
	\includegraphics{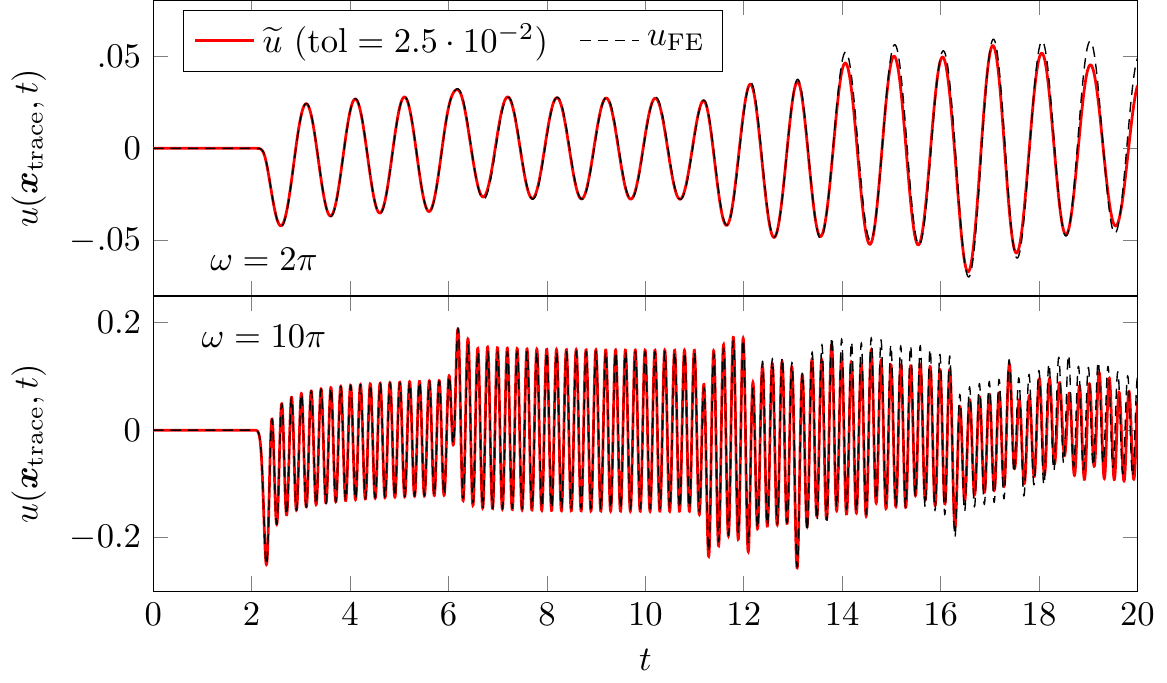}
	\caption{Value of solution at point $\vx_{\textup{trace}}=(-1,-2)$ for excitation frequencies $\omega=2\pi$ (top) and $\omega=10\pi$ (bottom). The reference FE solution is also included in both plots as a dashed line.}
	\label{fig:room:h:trace}
\end{figure}

\begin{figure}[tb]
	\centering
	\includegraphics[trim={0 0 0 2cm},clip]{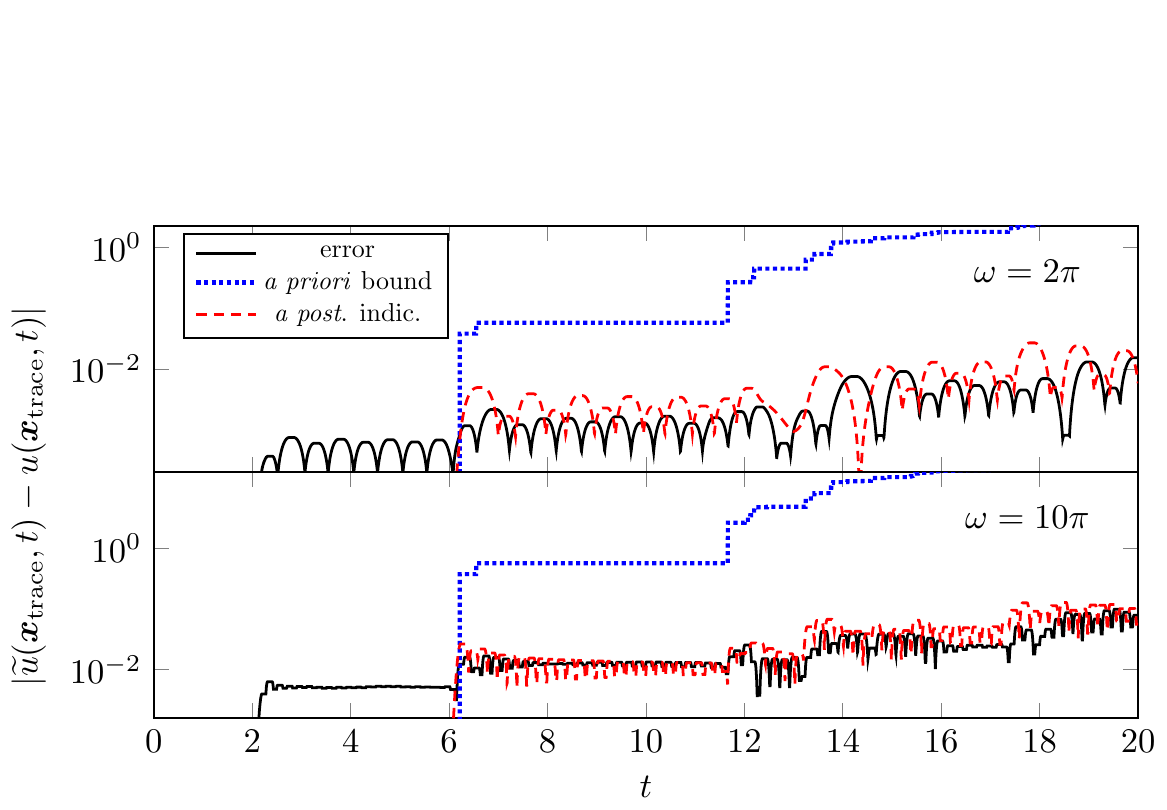}
	\caption{Approximation error at $\vx_{\textup{trace}}=(-1,-2)$ for excitation frequencies $\omega=2\pi$ (top) and $\omega=10\pi$ (bottom). We also include the \emph{a priori} error bound and the heuristic \emph{a posteriori} error indicator using dotted and dashed lines, respectively. Note that, to avoid clutter, we show a ``filtered'' version of the error, using a window of width $0.1$: $\max_{|\delta|\leq0.05}|\uapp(\vx_{\textup{trace}},t+\delta)-u(\vx_{\textup{trace}},t+\delta)|$.}
	\label{fig:room:h:err}
\end{figure}

We display the approximation error in \cref{fig:room:h:err}. There, we also display an \emph{a priori} error bound based on \cref{eq:err}, as well as the heuristic \emph{a posteriori} error indicator $|\uapp_{\text{GO}}-\uapp|$, with $\uapp_{\text{GO}}$ being the approximation obtained with our approach under the assumption of geometrical optics. The \emph{a priori} error bound is piecewise-constant since time-independent terms are added to the error whenever new diffraction waves reach $\vx_{\textup{trace}}$.

For $t<6$, no diffraction is present, the error indicators are both zero, since the only errors are due to discretization of the PDE\footnote{Discretization errors are invisible to the two indicators but can normally be estimated through other means, e.g., using standard results for FEM-based approximation of the wave equation.}. For $t>6$, the \emph{a priori} estimator provides a reliable but pessimistic upper bound for the error. On the other hand, the \emph{a posteriori} indicator manages to stay much closer to the actual error. Although not substantiated by theory, this error indicator seems quite effective at giving quantitative information on the error magnitude, without the need to compute a reference solution.

Note that using FE to approximate the wave $u$ requires carrying out a new simulation from scratch for every frequency to be studied. To this end, one needs to choose a mesh with $\omega$-dependent resolution: the mesh size should be small enough for the well-known \emph{pollution effect} (see, e.g., \cite{BabuskaSauter1997}) to be absent.

A constraint on the mesh resolution is present also in our proposed approach. However, it only applies to the problem defining the free-space solution $\Psi$, which is 1-dimensional in space. Hence, having to refine the mesh represents a much smaller obstacle to efficiency. In particular, as the frequency $\omega$ increases, the computation of $\uapp$ becomes comparatively more and more efficient, with respect to the computation of $\ufem$.

\section{Conclusions}\label{sec:conclusions}
We have presented a surrogate modeling strategy for approximating waves propagating through complex 2-dimensional domains with polygonal boundaries. Our method relies on the automatic identification of reflection and diffraction effects caused by the domain geometry. Each effect is modeled through a relatively simple nonlinear expression. Reflection-related components are built using geometrical optics, whereas diffraction-related components are modeled by a novel \emph{ad hoc} modification of the geometrical theory of diffraction. We have also provided estimators for the quality of the approximation achieved by our method.

In our numerical tests, we have observed a good accuracy, with the main features of the target wave being well identified. Notably, our diffraction model has proven to be fairly effective. Still, it relies on the parameter $\ks$, which, in some sense, determines the strength of surrogate diffraction effects. Although we have presented some heuristics for choosing $\ks$, more refined strategies for selecting $\ks$ and, more generally, a thorough validation of our diffraction model remain open issues.

In terms of complexity, our method requires the solution of a simplified 1D-in-space problem, much simpler than the original 2D-in-space one. Another favorable aspect of our algorithm is its potential to be run on parallel architectures, since the computation of different rays can be carried out independently. This is not the case for standard timestepping-based discretizations, due to their intrinsically sequential nature.

\subsection{Open questions and possible extensions}
We now present some open questions concerning our method, which give cues for further research on the subject.

\paragraph{Additional physics} A first question is whether our algorithm may be applied to wave equations incorporating further physical effects, like dissipation. In such setting, information can propagate faster than the wave speed due to diffusion. In our opinion, it should still be possible to apply our method to some effect, as long as dissipation is kept low, so as to maintain at least a partial notion of causality.

\paragraph{3D} Concerning other extensions, we have already mentioned that our geometry-based approach to wave propagation generalizes from 2D to 3D. In 3D, reflection and diffraction happen at facets and edges, as opposed to edges and vertices, respectively \cite{mcnamara1990introduction}. Moving from polygons to polyhedra leads to an increase in the ``optical entities'' to be considered (e.g., triangles have 3 edges and 3 vertices, whereas tetrahedra have 4 facets and 6 edges), especially if complex 2D surfaces need to be meshed. This leads to an increase in the number of field components to be considered. If the number of \emph{successive} reflections/diffractions is large, this effect is further amplified. For instance, if at most 3 (resp., 5) successive reflection or diffraction events happen involving tetrahedral scatterers, we may expect the cost of our method to increase by a factor of 5 (resp., 13) with respect to a similar simulation with triangular scatterers in 2D.

Moreover, some additional implementation issues arise in 3D due to the higher spatial dimension. For instance, when computing the spatial supports $\Omega_n$, our current strategy, cf.~\cref{rem:cone}, could be rather expensive if implemented naively, especially if the number of facets of $\partial\Omega$ is large. More sophisticated strategies to represent the sets $\Omega_n$, e.g., involving \emph{beam}-tracing \cite{svensson_savioja_overview} as opposed to \emph{ray}-tracing, may be required to handle 3 dimensions more effectively. That being said, besides the geometrical side, our method still only requires solving a 1D FEM problem in spherical symmetry, even though the physical dimension of the underlying problem has increased.

\paragraph{Coupling with other methods} Our method is not too well suited to approximate wave propagation in trapping media, since the number $N$ of terms in the approximation may increase dramatically. To alleviate this issue, one might consider splitting the domain into non-trapping and trapping parts. Our approach could be then applied to approximate wave propagation in non-trapping subdomains, while standard MOR approaches (based, e.g., on projection onto a basis of standing waves) may deliver an efficient approximation over the trapping region. As in other MOR methods based on domain decomposition \cite{DDROM1,DDROM2}, the most complex step is the (dynamic) coupling between subdomains. This remains an open issue for our method. We also note that, more generally, developing such an approach could also enable applying our method within multi-physics problems, e.g., when modeling the interaction between propagating wave and structures.

\paragraph{Parametrized problems} Finally, we recall that, in many applications, the ultimate target is understanding how the wave $u$ solving \cref{eq:wave} depends on underlying parameters $\mathbf{p}$, e.g., the forcing term $f$, the shape of the domain $\Omega$, etc. In this setting, MOR tries to construct a surrogate model of the form $\uapp=\uapp(\vx,t;\mathbf{p})$, providing a good approximation of $u$ over a whole range of parameter values. Even though our technique was presented here in the non-parametric setting, we believe that it potentially allows incorporating the parameter dependence in a natural and efficient way. In our opinion, this might be achievable by leveraging the simple and interpretable structure of the field components (free-space solution, spatial support, and angular modulation). As a simple preliminary example, we showcased this in \cref{sec:numerical:helm} for a parametric source term, with the parameter being the frequency. We are currently investigating how to extend our method to more complicated parametric problems.

\appendix

\section{UTD diffraction coefficient}\label{app:diffract}
\begin{figure}[tb]
	\centering
	\includegraphics{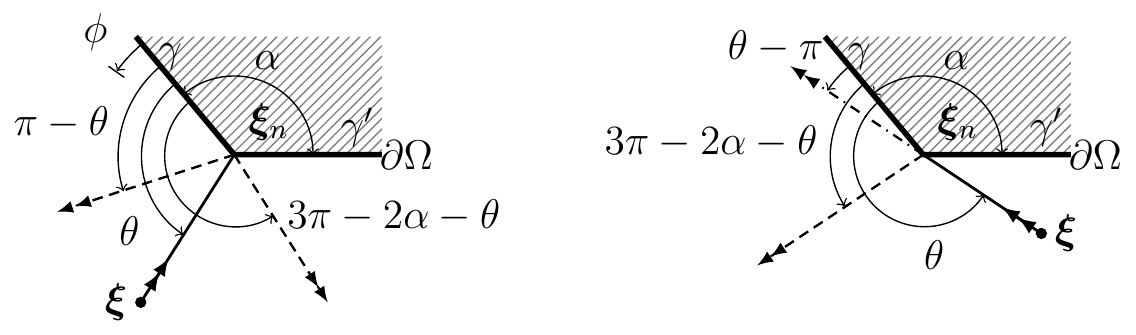}
	\caption{Diagrams for the two cases of scattering for concave corners ($0<\alpha<\pi$): without (left plot) and with shadow zone (right plot). The dashed lines are reflection shadow boundaries. The dash-dotted line is an incident shadow boundary. Shadow regions are absent if and only if $\pi-\alpha\leq\theta\leq\pi$. The angular coordinate $0<\phi<2\pi-\alpha$ is measured starting from one of the two adjacent edges of $\partial\Omega$.}
	\label{fig:scatter}
\end{figure}

Consider the setup shown in \cref{fig:scatter}: a wedge with exterior angle $\alpha$ and vertex $\vxi_n$ is hit by a wave coming from $\vxi$, located at an angular position $\phi=\theta$. The angle $\phi$ is measured from one of the two sides of the wedge. For simplicity, we assume that each of the wedge's sides is either sound-soft or sound-hard. For a discussion of the UTD coefficients in the case of impedance boundary conditions, we refer to \cite{Manara}.

Let $\vx$ be a point located at distance $s=\norm{\vx-\vxi_n}$ from $\vxi_n$, at an angle $\phi=\phi(\vx)$. The diffraction coefficient $D=D(\vx)$ 
represents 
the magnitude (and sign) of the diffraction wave at $\vx$, assuming that the incident wave, with origin at $\vxi$, is time-harmonic with unit amplitude and frequency $k$. UTD \cite{kouyoumjian1974uniform} predicts a diffraction coefficient
\begin{equation*}
    D=D_1+D_2\pm(D_3+D_4),
\end{equation*}
with the sign $\pm$ depending on the type of boundary conditions (``$+$'' for Neumann, ``$-$'' for Dirichlet). Given the wedge index $\nu:=\pi/(2\pi-\alpha)$, all contributions $D_j$ have the form
\begin{equation*}
    D_j= -\frac{\nu}{2\sqrt{2\pi ks}}\cot(\widehat\phi_j)F\left(2ks\cos(\widetilde\phi_j)^2\right),
\end{equation*}
with
\begin{align*}
    \widehat\phi_1=\frac{\pi+\phi-\theta}2\nu,\quad\widetilde\phi_1=&N_1(2\pi-\alpha)-\frac{\phi-\theta}2,\\
    \widehat\phi_2=\frac{\pi-\phi+\theta}2\nu,\quad\widetilde\phi_2=&N_2(2\pi-\alpha)-\frac{\phi-\theta}2,\\
    \widehat\phi_3=\frac{\pi+\phi+\theta}2\nu,\quad\widetilde\phi_3=&N_3(2\pi-\alpha)-\frac{\phi+\theta}2,\\
    \widehat\phi_4=\frac{\pi-\phi-\theta}2\nu,\quad\widetilde\phi_4=&N_4(2\pi-\alpha)-\frac{\phi+\theta}2,
\end{align*}
and
\begin{equation*}
    F(x)=2\sqrt{x}\sqrt{\left(\int_{\sqrt{x}}^\infty\cos(y^2)\textrm{d}y\right)^2+\left(\int_{\sqrt{x}}^\infty\sin(y^2)\textrm{d}y\right)^2}.
\end{equation*}
The integers $N_j$ are defined as
\begin{equation*}
    N_1=\left[\frac\nu2\right],\quad N_2=\left[-\frac\nu2\right],\quad N_3=\left[\frac{1+\nu}2\right],\quad\text{and}\quad N_4=\left[\frac{1-\nu}2\right],
\end{equation*}
with square brackets denoting the rounding operator, which yields the integer that is closest to its argument.

The diffraction coefficient has unit jumps (with suitable signs) at the \emph{shadow boundaries}, i.e., at the angular positions corresponding to boundaries of the spatial support of the source wave (in which case, we have an ``incident shadow boundary'', ISB) or of one of its reflections (in which case, we have a ``reflection shadow boundary'', RSB). The locations of the shadow boundaries can be identified geometrically. For instance, for convex wedges ($0<\alpha<\pi$), shadow boundaries are at $\phi=\abs{\pi-\theta}$ (which is an ISB if $\theta>\pi$, and an RSB otherwise) and at $\phi=2\pi-\alpha-\abs{\pi-\alpha-\theta}$ (which is an ISB if $\theta<\pi-\alpha$, and an RSB otherwise).

\begin{figure}[tb]
	\centering
	\includegraphics{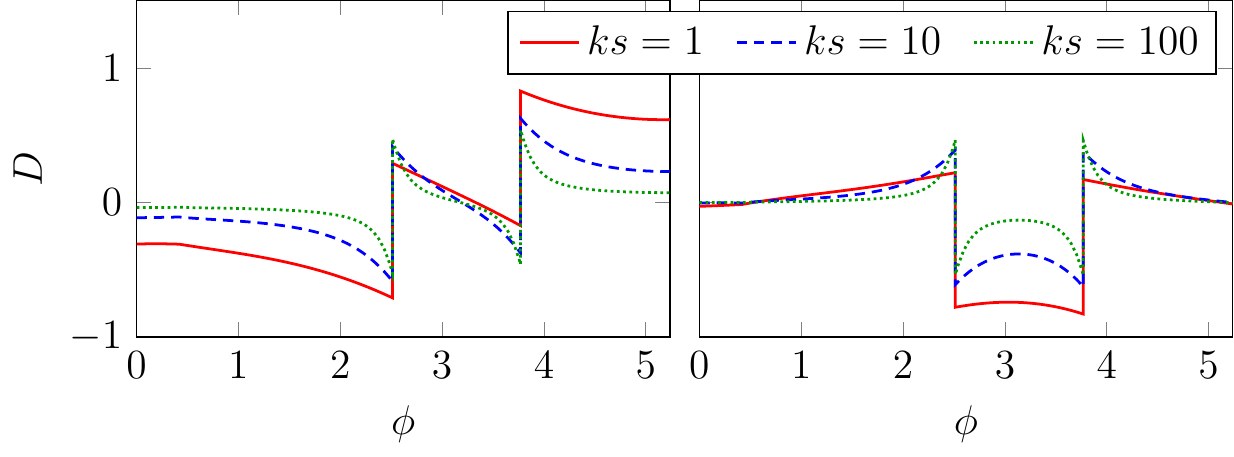}
	\caption{UTD diffraction coefficients for Neumann (left) and Dirichlet (right) wedges.}
	\label{fig:diff_coeff}
\end{figure}

The diffraction coefficient for $\alpha=\frac\pi3$, $\theta=\frac\pi5$, Neumann conditions, and $ks\in\{1,10,100\}$ are shown in \cref{fig:diff_coeff} (left). Two jumps happen, at an RSB (at $\phi=\frac45\pi$) and at an ISB (at $\phi=\frac65\pi$). Similar results are obtained for Dirichlet conditions, cf.~\cref{fig:diff_coeff} (right).

We can observe that the value of $ks$ determines how quickly $D$ decays to 0 around each shadow boundaries: larger values of $ks$ yield a narrower ``transition region'', which, loosely speaking, is the support of $D$ around each discontinuity. Indeed, asymptotically in $ks$, the width of each transition region is proportional to $1/\sqrt{ks}$ \cite{mcnamara1990introduction}.

\bibliographystyle{elsarticle-num} 
\bibliography{references.bib}
\end{document}